\documentclass{amsart}
\newtheorem{theorem}{Theorem}[subsection]
\newtheorem{lemma}[theorem]{Lemma}
\newtheorem{proposition}[theorem]{Proposition}
\newtheorem{corollary}[theorem]{Corollary}
\theoremstyle{definition}
\newtheorem{definition}[theorem]{Definition}
\newtheorem{example}[theorem]{Example}

\theoremstyle{remark}
\newtheorem{remark}[theorem]{Remark}

\newcommand{\Tr}{\mbox{Tr}\,} 
\newcommand{\tr}{\mbox{tr}}
\newcommand{\id}{\mbox{id}} 
 
\newcommand{\Ker}{\mbox{Ker\,}} 
 
\newcommand{\Index}{\mbox{Index}\,} 
\newcommand{\End}{\mbox{End}} 
\newcommand{\URep}{\mbox{URep}}
\newcommand{\Rep}{\mbox{Rep}}
\newcommand{\Corep}{\mbox{Corep}}
\newcommand{\Bimod}{\mbox{Bimod}}
\newcommand{\Hom}{\mbox{Hom}} 
\newcommand{\Map}{\mbox{Map}}
\newcommand{\Coinv}{\mbox{Coinv}}
 
\renewcommand{\span}{\mbox{span}} 
    
\newcommand{\Ad}{\mbox{Ad}\,}
\newcommand{\R}{{\mathcal R}}
\newcommand{\bR}{{\bar \R}}
\newcommand{\bTheta}{{\bar \Theta}}
\newcommand{\g}{\mathfrak{g}}
\newcommand{\bF}{{\bar F}}

\newcommand{\tJ}{{\mathcal J}}
\newcommand{\C}{{\mathcal C}}

\newcommand{\eps}{\varepsilon} 
\newcommand{\op}{op}
\newcommand{\cop}{cop}

\newcommand{\act}{\rightharpoonup} 
\newcommand{\lact}{\triangleright} 
\newcommand{\actr}{\rightharpoonup} 
\newcommand{\actl}{\leftharpoonup}

\newcommand{\ract}{\triangleleft} 
\newcommand{\la}{\langle\,} 
\newcommand{\ra}{\,\rangle}

\renewcommand{\star}{{\dag}} 
\newcommand{\rtimes}{{>\!\!\!\triangleleft}} 
\newcommand{\ltimes}{{\triangleright\!\!\!<}} 
\newcommand{\smrtimes}{>\!\!\triangleleft} 
 
\newcommand{\half}{{1/2}}
 
\newcommand{\0}{_{(0)}}  
\newcommand{\1}{_{(1)}} 
\newcommand{\2}{_{(2)}} 
\newcommand{\3}{_{(3)}} 
\newcommand{\4}{_{(4)}} 
\renewcommand{\0}{^{(0)}} 
\newcommand{\I}{^{(1)}} 
\newcommand{\II}{^{(2)}}   
\newcommand{\III}{^{(3)}}

\begin{document}

\title{\bf Finite Quantum Groupoids and Their Applications}

\author{Dmitri Nikshych}
\address{UCLA, Department of Mathematics, 405 Hilgard Avenue,
Los Angeles, CA 90095-1555, USA}
\email{nikshych@math.ucla.edu}
\thanks{The first author thanks P.~Etingof for numerous stimulating
discussions on quantum groupoids and to MIT for the hospitality 
during his visit}

\author{Leonid Vainerman}
\address{Universit\'e de Strasbourg, D\'epartement de
Math\'ematiques, 7, rue Ren\'e Descartes, F-67084 Strasbourg, France}
\email{wain@agrosys.kiev.ua, vaynerma@math.u-strasbg.fr}
\thanks{ The second author is grateful to M.~Enock, 
V.~Turaev and J.-M.~Vallin for many valuable discussions and to 
l'Universit\'e Louis Pasteur (Strasbourg) 
for the kind hospitality during his work on this survey}

\date{June 7, 2000}

\begin{abstract}
We give a survey of the theory of finite quantum groupoids 
(weak Hopf algebras), including foundations of the theory and applications 
to finite depth  subfactors,  dynamical deformations of quantum groups,  
and invariants of knots and 3-manifolds.
\end{abstract}

\maketitle
\tableofcontents

\begin{section}
{Introduction}
By {\em quantum groupoids} we understand
weak Hopf algebras introduced in \cite{BNSz}, \cite{BSz1}, \cite{N}.
These objects generalize Hopf algebras (in fact, many Hopf-algebraic
concepts can be extended to the quantum groupoid case) and usual
finite groupoids (see \cite{NV1} for discussion). 
Every quantum groupoid has two canonical subalgebras
that play the same role as the space of units in a groupoid
and projections on these
subalgebras generalizing the source and target maps in a groupoid.
We use the term ``quantum groupoid'' instead of
``weak Hopf algebra'' to stress this similarity that leads to many 
interesting constructions and examples.

Our initial motivation for studying quantum groupoids in
\cite{NV1}, \cite{NV2},\cite{N1} was their connection
with depth 2 von Neumann subfactors, first mentioned in \cite{O2},
which was also one of the main topics of \cite{BNSz}, \cite{BSz1}, 
\cite{BSz2}, \cite{NSzW}. It was shown in \cite{NV2} that 
quantum groupoids naturally arise as non-commutative symmetries
of subfactors, namely if $N\subset M\subset M_1\subset M_2\subset\dots$
is the Jones tower constructed from a finite index depth $2$ 
inclusion $N\subset M$ of II$_1$ factors, then $H=M'\cap M_2$ 
has a canonical structure of a quantum groupoid
acting outerly on $M_1$ such that $M =M_1^H$ and
$M_2 = M_1\rtimes H$, moreover $\widehat H=N'\cap M_1$ is a
quantum groupoid dual to $H$.

In \cite{NV3} we extended this result to show that
quantum groupoids give a description of arbitrary
finite index and finite depth II$_1$ subfactors via
a Galois correspondence and
explained how to express subfactor invariants such as bimodule   
categories and principal graphs in quantum groupoid terms. Thus,
in this respect quantum groupoids play the same role as Ocneanu's paragroups 
\cite{O1}, \cite{D}.    

Quantum groupoids also appear naturally in the theory of
dynamical deformations of quantum groups \cite{EV}, \cite{Xu}.
It was shown in \cite{EN} that for every simple Lie algebra
$\g$ the corresponding Drinfeld-Jimbo quantum group $U_q(\g)$
for $q$ a primitive root of unity has a family of dynamical
twists that give rise to self-dual finite
quantum groupoids. To every such a twist one can associate
a solution of the quantum dynamical Yang-Baxter equation \cite{ES}.

As in the case of Hopf algebras, the representation category $\Rep(H)$
of a quantum groupoid $H$ is a monoidal category with duality. Existence
of additional (quasitriangular, ribbon, or modular) structures on $H$ 
makes $\Rep(H)$, respectively,  braided, ribbon, or modular \cite{NTV}.
It is well known \cite{RT}, \cite{T}, that such categories give rise to
invariants of knots and 3-manifolds. 

The survey is organized as follows. 

In Section 2 we give basic definitions and examples of finite 
quantum groupoids, discuss their fundamental properties 
(following \cite{BNSz} and \cite{NV1}) and relations with other versions
of quantum groupoids. 

The exposition of the theory of integrals in
Section 3 follows \cite{BNSz} and is similar to the Hopf algebra case 
\cite{M}. The main results here are extensions of the fundamental 
theorem for Hopf modules and  the  Maschke theorem. The latter gives 
an equivalence between existence of normalized integrals and 
semi-simplicity of a quantum groupoid.

In  Section 4 we define the notions of action and smash product
and extend the well known Blattner-Montgomery
duality theorem for Hopf algebra actions to quantum groupoids 
(see \cite{N2}). 

In Section 5 we describe monoidal structure and duality on the representation 
category of a quantum groupoid (\cite{BSz2} and \cite{NTV}). Then we 
distinguish special cases of quantum groupoids leading to
braided, ribbon, and modular categories of representations. 

We also develop the Drinfeld
double construction and emphasize the role of factorizable quantum
groupoids in producing modular categories (see \cite{NTV}).
Section 6  contains an extension of the Drinfeld twisting procedure to 
quantum groupoids and construction of a quantum groupoid via  a
{\em dynamical twisting} of a Hopf algebra. Important concrete examples
of such twistings are {\em dynamical quantum groups} at roots of 1 \cite{EN}. 

Semisimple and $C^*$-quantum groupoids are analyzed in
Section 7. Following \cite{BNSz}, our discussion is concentrated 
around the existence of the Haar integral in connection with a
special implementation of the square of the antipode.

Sections 8 and 9 are devoted to the description of finite index an
finite depth subfactors via crossed products of II$_1$ factors
with $C^*$-quantum groupoids. These results  
were obtained in \cite{NV2}, \cite{NV3} and partially in 
\cite{N1}, \cite{NSzW}.

\end{section}

\begin{section}
{Definitions and Examples}

In this section we give definitions and discuss basic
properties of finite quantum groupoids. Most of the material 
presented here can be found in \cite{BNSz} and \cite{NV1}. 

Throughout this paper we use Sweedler's notation for
comultiplication, writing $\Delta(b) = b\1 \otimes b\2$.
Let $k$ be a field.

\subsection{Definition of the quantum groupoid}
By {\em quantum groupoids} or {\em weak Hopf algebras}
we understand the objects introduced in \cite{BNSz}, \cite{BSz1}
as a generalization of ordinary Hopf algebras. 

\begin{definition} 
\label{finite quantum groupoid}
A {\em (finite) quantum groupoid} over $k$ is a
finite dimensional $k$-vector space $H$
with the structures of an associative algebra $(H,\,m,\,1)$ 
with multiplication $m:H\otimes_k H\to H$ and unit $1\in H$ and a 
coassociative coalgebra $(H,\,\Delta,\,\eps)$ with comultiplication
$\Delta:H\to H\otimes_k H$ and counit $\eps:H\to k$ such that:
\begin{enumerate}
\item[(i)] The comultiplication $\Delta$ 
is a (not necessarily unit-preserving) homomorphism of algebras such that
\begin{equation}
(\Delta \otimes \id) \Delta(1) =
(\Delta(1)\otimes 1)(1\otimes \Delta(1)) =
(1\otimes \Delta(1))(\Delta(1)\otimes 1),
\label{Delta 1}
\end{equation}
\item[(ii)] The counit is a $k$-linear map
satisfying the identity:
\begin{equation}
\eps(fgh) = \eps(fg\1)\, \eps(g\2h) = \eps(fg\2)\, \eps(g\1h),
\label{eps m}
\end{equation}
for all $f,g,h\in H$. 
\item[(iii)]
There is a linear map $S: H \to H$, called an {\em antipode}, such that,
for all $h\in H$,
\begin{eqnarray}
m(\id \otimes S)\Delta(h) &=&(\eps\otimes\id)(\Delta(1)(h\otimes 1)),
\label{S epst} \\
m(S\otimes \id)\Delta(h) &=& (\id \otimes \eps)((1\otimes h)\Delta(1)),
\label{S epss} \\
m(m\otimes id)(S\otimes \id \otimes S)(\Delta\otimes id)\Delta(h) &=& S(h).
\label{S id S}
\end{eqnarray}
\end{enumerate} 
\end{definition}
A quantum groupoid is a Hopf algebra if and only if the comultiplication 
is unit-preserving if and only if the counit is a homomorphism of algebras.

A {\em morphism} between quantum groupoids $H_1$ and $H_2$
is a map $\alpha : H_1 \to H_2$ which is both algebra and coalgebra
homomorphism preserving unit and
counit and which intertwines the antipodes of $H_1$ and $H_2$,
i.e., $ \alpha\circ S_1 = S_2\circ \alpha$. The image of a morphism is
clearly a quantum groupoid. The tensor product of two quantum groupoids
is defined in an obvious way. 

\subsection{Counital maps and subalgebras}
The linear maps defined in (\ref{S epst}) and (\ref{S epss})
are called {\em target} and {\em source counital maps}
(see examples below for explanation of the terminology) and 
denoted $\eps_t$ and $\eps_s$ respectively :
\begin{equation}
\eps_t(h) = (\eps\otimes\id)(\Delta(1)(h\otimes 1)),\qquad
\eps_s(h) = (\id \otimes \eps)((1\otimes h)\Delta(1)).
\end{equation}

In the next proposition we collect several useful properties of 
the counital maps.

\begin{proposition}
\label{counital properties}
For all $h,g \in H$ we have 
\begin{enumerate}
\item[(i)] Counital maps are idempotents in $\End_k(H)$:   
$$
\eps_t(\eps_t(h))=\eps_t(h), \qquad  \eps_s(\eps_s(h))=\eps_s(h),
$$
\item[(ii)] the relations between $\eps_t,\,\eps_s$, and comultiplication
are as follows
$$
(\id \otimes \eps_t)\Delta(h) =1\1h\otimes 1\2, \qquad
(\eps_s \otimes \id)\Delta(h) =1\1 \otimes h1\2,
$$
\item[(iii)] the images of counital maps are characterized by
$$
h = \eps_t(h) \mbox{ iff } \Delta(h) = 1\1h\otimes 1\2, \qquad
h = \eps_s(h) \mbox{ iff }  \Delta(h) =1\1 \otimes h1\2,
$$
\item[(iv)]  $\eps_t(H)$ and $\eps_s(H)$  commute,
\item[(v)] one also has identities dual to (ii) :
$$
h\eps_t(g) = \eps(h\1g)h\2, \qquad \eps_s(h)g = h\1\eps(gh\2).
$$
\end{enumerate}
\end{proposition}
\begin{proof}
We prove the identities containing the target counital map, 
the proofs of their source counterparts are similar. Using the axioms 
(\ref{Delta 1}) and (\ref{eps m}) we compute 
$$
\eps_t(\eps_t(h))=\eps(1\1h)\eps(1'\1 1\2)1'\2 =
\eps(1\1h)\eps(1\2) 1\3 = \eps_t(h),
$$ where $1'$ stands for the second copy of the unit, proving (i).
For (ii) we have
\begin{eqnarray*}
h\1 \otimes \eps_t(h\2)
&=& h\1 \eps(1\1h\2) \otimes 1\2 \\
&=& 1\1h\1 \eps(1\2h\2) \otimes 1\3 = 1\1h \otimes 1\2.  
\end{eqnarray*}
To prove (iii) we observe that
$$
\Delta(\eps_t(h)) = \eps(1\1h)1\2 \otimes 1\3
= \eps(1\1h) 1'\1 1\2 \otimes 1'\2 = 1'\1  \eps_t(h) \otimes 1'\2, 
$$
on the other hand, applying $(\eps\otimes \id)$ to both sides of
$\Delta(h) =1\1h \otimes 1\2$, we get $h=\eps_t(h)$. (iv) is immediate
in view of the identity $1\1 \otimes 1'\1 1\2 \otimes 1'\2 =
1\1 \otimes 1\2  1'\1 \otimes 1'\2 $. Finally, we show (v) :
\begin{eqnarray*}
\eps(h\1g)h\2 
&=& \eps(h\1g\1)h\2 \eps_t(g\2) = \eps(h\1g\1)h\2 g\2 S(g\3) \\
&=& hg\1 S(g\2) = h\eps_t(g),
\end{eqnarray*}
where the antipode axiom (\ref{S epst}) is used.
\end{proof}

The images of the counital maps
\begin{eqnarray}
H_t &=& \eps_t(H)  =  \{h\in H \mid \Delta(h) = 1\1 h\otimes 1\2 \}, \\
H_s &=& \eps_s(H)  = \{h\in H \mid \Delta(h) = 1\1 \otimes h 1\2 \}
\end{eqnarray}
play the role of bases of $H$. The next proposition
summarizes their properties.

\begin{proposition}
\label{base properties}
$H_t$ (resp.\ $H_s$) is a left (resp.\ right) coideal subalgebra
of $H$. These subalgebras commute with each other; moreover
\begin{equation*}
H_t = \{(\phi\otimes \id)\Delta(1) \mid \phi\in \widehat H \}, \qquad
H_s = \{(\id \otimes \phi)\Delta(1) \mid \phi\in \widehat H \},
\end{equation*}
i.e., $H_t$ (resp.\ $H_s$) is generated by the right (resp.\ left)
tensorands of $\Delta(1)$. 
\end{proposition}
\begin{proof}
$H_t$ and $H_s$ are coideals
by  Proposition~\ref{counital properties}(iii),
they commute by \ref{counital properties}(iv). 
We have   $H_t =\eps(1\1 H)1\2 \subset 
\{(\phi\otimes \id)\Delta(1) \mid \phi\in \widehat H \}$, conversely
$\phi(1\1)1\2 = \phi(1\1)\eps_t(1\2) \subset H_t$, therefore
$H_t = \{(\phi\otimes \id)\Delta(1) \mid \phi\in \widehat H \}$.
To see that it is an algebra we note that $1=\eps_t(1) \in H_t$
and for all $h,g\in H$ compute, using 
Proposition~\ref{counital properties}(ii) and (v) :
\begin{eqnarray*}
\eps_t(h)\eps_t(g) 
&=& \eps (\eps_t(h)\1 g) \eps_t(h)\2 \\
&=& \eps (1\1 \eps_t(h)g)1\2 = \eps_t(\eps_t(h)g) \in H_t.
\end{eqnarray*}
The statements about $H_s$ are proven similarly.
\end{proof}

\begin{definition}
\label{counital subalgebras}
We call $H_t$ (resp.\ $H_s$) a {\em target} (resp.\ {\em source})
{\em  counital subalgebra}. 
\end{definition}

\subsection{Properties of the antipode}
The properties of the antipode of a quantum groupoid 
are similar to those of a finite-dimensional Hopf algebra.

\begin{proposition}[\cite{BNSz}, 2.10]
\label{properties of S}
The antipode $S$ is unique and bijective. Also,
it is both algebra and coalgebra anti-homomorphism. 
\end{proposition}
\begin{proof}
Let $f*g = m(f\otimes g)\Delta$ be the convolution of $f,g \in 
\End_k(H)$. Then $S*\id =\eps_s,\, \id*S=\eps_t$, and $S*\id*S=S$. 
If $S'$ is another antipode of $H$ then
$$
S' = S' *\id * S' = S' *\id * S = S*\id*S=S.
$$
To check that $S$ is an algebra anti-homomorphism, we compute
\begin{eqnarray*}
S(1) 
&=& S(1\1)1\2S(1\3)= S(1\1)\eps_t(1\2) =\eps_t(1) =1, \\
S(hg) 
&=& S(h\1g\1)\eps_t(h\2g\2) = S(h\1g\1) h\2 \eps_t(g\2) S(h\3) \\
&=& \eps_s(h\1g\1)S(g\2)S(h\2) = S(g\1) \eps_s(h\1) \eps_t(g\2) S(h\2) 
   = S(g)S(h),
\end{eqnarray*}
for all $h,g\in H$, where we used Proposition~\ref{counital properties}(iv)
and easy identities $ \eps_t(hg) = \eps_t(h\eps_t(g))$ and
$ \eps_s(hg) = \eps_t(\eps_s(h)g)$.
Dualizing the above arguments we show that $S$ is also
a coalgebra anti-homomorphism :
\begin{eqnarray*}
\eps(S(h))
&=& \eps(S(h\1)\eps_t(h\2)) = \eps (S(h\1)h\2) =\eps(\eps_t(h)) =\eps(h),\\
\Delta(S(h))
&=& \Delta(S(h\1)\eps_t(h\2))\\
&=& \Delta(S(h\1)) (\eps_t(h\2) \otimes 1) \\
&=& \Delta(S(h\1)) (h\2 S(h\4) \otimes \eps_t(h\3)) \\
&=& \Delta(\eps_s(h\1)) ( S(h\3) \otimes S(h\2) ) \\
&=& S(h\3) \otimes \eps_s(h\1)S(h\2) = S(h\2) \otimes S(h\1). 
\end{eqnarray*}
The proof of the bijectivity of $S$ can be found in (\cite{BNSz}, 2.10).
\end{proof}

Next, we investigate the relations between the antipode
and counital maps.

\begin{proposition}
\label{S and counital maps}
We have $S\circ \eps_s = \eps_t\circ S$ and 
$ \eps_s\circ S = S\circ \eps_t$. The restriction of $S$
defines an algebra anti-isomorphism between counital
subalgebras $H_t$ and $H_s$. 
\end{proposition}
\begin{proof}
Using results of Proposition~\ref{properties of S}
we compute
$$
S(\eps_s(h)) = S(1\1) \eps(h1\2) = \eps(1\1S(h))1\2 =\eps_t(S(h)),
$$
for all $h\in H$. The second identity is proven similarly. Clearly,
$S$ maps $H_t$ to $H_s$ and vice versa. Since $S$ is bijective,
and $\dim H_t =\dim H_s$ by Proposition~\ref{base properties},
therefore $S|_{H_t}$ and  $S|_{H_s}$ are anti-isomorphisms. 
\end{proof}

\begin{proposition}[\cite{NV1}, 2.1.12]
\label{full subcategory}
Any nonzero morphism $\alpha :H \to K$ of quantum groupoids preserves
counital subalgebras, i.e. $H_t \cong K_t$ and $H_s \cong K_s$.
Thus, quantum groupoids with a given target (source) counital
subalgebra form a full subcategory.
\end{proposition}
\begin{proof}
It is clear that $\alpha|_{H_t} : H_t \to K_t$ is a homomorphism. 
If we write 
$$
\Delta(1_H) = \sum_{i=1}^n\, w_i\otimes z_i
$$
with $\{ w_i\}_{i=1}^n$ and $\{ z_i\}_{i=1}^n$ linearly independent,
then $\Delta(1_K) = \sum_{i=1}^n\, \alpha(w_i)$ $\otimes \alpha(z_i)$.
By Proposition~\ref{base properties}, $K_t =\span\{\alpha(z_i)\}$,
i.e., $\alpha|_{H_t}$ is surjective. 
Since 
$$
z_j = \eps_t(z_j) = \sum_{i=1}^n\, \eps(w_iz_j)z_i,
$$ 
then  $\eps(w_iz_j) =\delta_{ij}$, therefore, 
\begin{eqnarray*}
\dim H_t  
&=& n = \sum_{i=1}^n\,\eps_H(w_i z_i) =  \sum_{i=1}^n\,\eps_H(w_i S(z_i)) \\
&=& \eps_H(\eps_t(1_H)) = \eps_H(1_H) = \eps_K(1_K) =  \dim K_t,
\end{eqnarray*}
so $\alpha|_{H_t}$ is bijective. The proof for source subalgebras is 
similar. 
\end{proof}

Let us recall that a
$k$-algebra $A$ is said to be {\em separable} \cite{P}
if the multiplication epimorphism $m: A\otimes_k A \to A$
has a right inverse as an $A-A$ bimodule homomorphism.
This is equivalent to the existence of a {\em separability element}
$e\in A\otimes_k A $ such that  $m(e) =1$ and 
$(a\otimes 1)e = e(1\otimes a),\ (1\otimes a)e=e(a\otimes 1)$
for all  $a\in A$.

\begin{proposition}
\label{separability}
The counital subalgebras $H_t$ and $H_s$ are separable,
their separability elements are
$e_t = (S \otimes \id)\Delta(1)$ and $e_s = (\id \otimes S)\Delta(1)$,
respectively.
\end{proposition}
\begin{proof}
For all $z\in H_t$ we compute, using Propositions~\ref{counital properties}
and \ref{S and counital maps} :
\begin{eqnarray*}
1\1 S^{-1}(z) \otimes 1\2
&=& S^{-1}(z)\1 \otimes \eps_t(S^{-1}(z)\2) \\
&=& 1\1 \otimes \eps_t(1\2S^{-1}(z)) 
    = 1\1 \otimes \eps_t(1\2z) = 1\1 \otimes 1\2z,
\end{eqnarray*}
applying $(S\otimes \id)$ to this identity we get $e_tz =ze_t$.
Clearly, $m(e_t)=1$, whence $e_t$ is a separability element.
The second statement follows similarly.
\end{proof}

\subsection{The dual quantum groupoid}
The set of axioms of Definition~\ref{finite quantum groupoid} is self-dual.
This allows to define a natural quantum groupoid
structure on the dual vector space $\widehat H=\Hom_k(H,k)$ by
``reversing the arrows'':
\begin{eqnarray}
& & \la h,\,\phi\psi \ra = \la \Delta(h),\,\phi\otimes\psi \ra, \\
& & \la h\otimes g,\,{\widehat\Delta}(\phi) \ra =  \la hg,\, \phi\ra, \\
& & \la h,\, {\widehat S}(\phi) \ra = \la S(h),\,\phi \ra,
\end{eqnarray}
for all $\phi,\psi \in \widehat H,\, h,g\in H$. The unit $\widehat 1$
of $\widehat H$ is $\eps$  and counit $\widehat\eps$ is
$\phi \mapsto \la\phi,\, 1\ra$.

In what follows we will use the Sweedler arrows,
writing for all $h\in H,\phi\in \widehat H$ :
\begin{equation}
h\actr\phi = \phi\1 \la h,\, \phi\2\ra,
\qquad
\phi\actl h =\la h,\,\phi\1 \ra \phi\2.
\end{equation}   
for all $h\in H,\phi\in \widehat H$.

The counital subalgebras of $\widehat H$ are canonically anti-isomorphic 
to those of $H$. More precisely, the map 
$H_t \ni z \mapsto (z\act \eps) \in {\widehat H}_s$ is an algebra isomorphism 
with the inverse given by  $\chi \mapsto (1 \actl \chi)$.
Similarly, the map $H_s \ni z \mapsto (\eps\actl z) \in {\widehat H}_t$
is an algebra isomorphism (\cite{BNSz}, 2.6).

\begin{remark}
\label{opposite groupoid}
The opposite algebra $H^{op}$ is also a quantum groupoid with the 
same coalgebra structure and the antipode $S^{-1}$. Indeed,
\begin{eqnarray*}
S^{-1}(h\2)h\1
&=& S^{-1}(\eps_s(h)) = S^{-1}(1\1) \eps(h1\2) \\
&=& S^{-1}(1\1) \eps(h S^{-1}(1\2)) =\eps(h1\1)1\2, \\
h\2 S^{-1}(h\1)
&=& S^{-1}(\eps_t(h)) = \eps(1\1h) S^{-1}(1\2) \\
&=& \eps( S^{-1}(1\1)h)S^{-1}(1\2) = 1\1\eps(1\2h),\\
S^{-1}(h\3) h\2 S^{-1}(h\1) 
&=& S^{-1}(h\1 S(h\2) h\3) = S^{-1}(h). 
\end{eqnarray*}
Similarly, the co-opposite coalgebra
$H^{cop}$ (with the same algebra structure 
as $H$ and the antipode $S^{-1}$) and $(H^{op/cop}, S)$ are
quantum groupoids.
\end{remark}

\subsection{Examples : groupoid algebras and their duals}
\label{groupoid algebras}

As group algebras and their duals are the easiest examples of Hopf
algebras, groupoid algebras and their duals provide 
examples of quantum groupoids (\cite{NV1}, 2.1.4).

Let $G$ be a finite {\em groupoid} (a category with finitely many
morphisms, such that each morphism is invertible), then the groupoid
algebra $kG$ (generated by morphisms $g\in G$ with the product of 
two morphisms being equal to  their composition
if the latter is defined and $0$ otherwise) 
is a quantum groupoid via :
\begin{equation}
\Delta(g) = g\otimes g,\quad \eps(g) =1,\quad S(g)=g^{-1},\quad g\in G.
\label{groupoid algebra}
\end{equation}
The counital subalgebras of $kG$ are equal to each other and
coincide with the abelian algebra spanned by the identity morphisms :
$(kG)_t = (kG)_s = \span\{gg^{-1}\mid g\in G\}$. The target and source
counital maps are given by the operations of taking the target (resp.\
source) object of a morphism : 
\begin{equation*}
\eps_t(g) =gg^{-1} = \id_{target(g)} \quad \mbox{ and } \quad
\eps_s(g) = g^{-1}g = \id_{source(g)}.
\end{equation*} 

The dual quantum groupoid $\widehat{kG}$ is
isomorphic to the algebra of functions on $G$, i.e.,
it is generated by idempotents  $p_g,\, g\in G$ such that
$p_g p_h= \delta_{g,h}p_g$, with the following structure operations
\begin{equation}
\Delta(p_g) =\sum_{uv=g}\,p_u\otimes p_v,\quad \eps(p_g)= \delta_{g,gg^{-1}},
\quad S(p_g) =p_{g^{-1}}.
\label{dual groupoid algebra}
\end{equation}
The target (resp.\ source) counital subalgebra is precisely the algebra
of functions constant on each set of morphisms of $G$ having
the same target (resp.\ source) object. The target and source maps are
\begin{equation*}
\eps_t(p_g) = \sum_{vv^{-1}=g}\, p_v  \quad \mbox{ and }  \quad  
\eps_s(p_g) = \sum_{v^{-1}v=g}\, p_v.
\end{equation*}

\subsection{Examples : quantum transformation groupoids}

It is known that any group action on a set (i.e., on a commutative algebra 
of functions) gives rise to a groupoid \cite{R}. 
Extending this construction, we associate a quantum groupoid 
with any action of a Hopf algebra on a separable algebra 
(``finite quantum space'').

Namely, let $H$ be a Hopf algebra and $B$ be a separable 
(and, therefore, finite dimensional and  semisimple \cite{P}) 
right $H$-module algebra with the action $ b\otimes h \mapsto b\cdot h$,
where $b\in B,h\in H$. Then $B^{op}$, the algebra opposite to $B$,
becomes a left $H$-module algebra via $h\otimes a \mapsto 
h\cdot a = a\cdot S_H(h)$.
One can form a {\em double crossed product algebra} $B^{op}\rtimes H\ltimes B$
on the vector space $ B^{op}\otimes H \otimes B$ with the multiplication
\begin{equation*}
(a\otimes h\otimes b)(a'\otimes h'\otimes b') =
(h\1\cdot a')a \otimes h\2{h'}\1 \otimes (b\cdot {h'}\2) b',
\end{equation*}
for all $a,a'\in B^{op},\, b,b'\in B,$ and $h,h'\in H$.

Assume that $k$ is algebraically closed and let
$e$ be the separability element of $B$ 
(note that $e$  is an idempotent when considered in $B\otimes B^{op}$).
Let $\omega\in B^*$ be uniquely determined by 
$(\omega\otimes \id)e = (\id\otimes \omega)e = 1$.

One can check that $\omega$ is the trace of the left
regular representation of $B$ and verifies 
the following identities:
$$
\omega((h\cdot a)b) = \omega(a (b\cdot h)), \qquad
e\I \otimes (h\cdot e\II)  = (e\I \cdot h) \otimes e\II,
$$
where $a\in B^{op}, b\in B$, and $e =e\I\otimes e\II$.

The structure of a quantum groupoid on  $B^{op}\rtimes H\ltimes B$
is given by
\begin{eqnarray} 
\Delta(a\otimes h\otimes b) &=& 
(a\otimes h\1 \otimes e\I)\otimes ((h\2\cdot e\II)\otimes h\3 \otimes b),\\
\eps(a\otimes h\otimes b) &=& \omega(a(h\cdot b)) = \omega(a(b\cdot S(h))),\\
S(a\otimes h\otimes b) &=& b\otimes S(h)\otimes a.
\end{eqnarray}  

\subsection{Examples : Temperley-Lieb algebras}
It was shown in \cite{NV3} (see Section $9$ for details) that any inclusion 
of type II${}_1$ factors with finite index and depth (\cite{GHJ}, 4.1) 
gives rise  to a  quantum groupoid describing the symmetry
of this inclusion. In the case of {\em Temperley-Lieb algebras} 
(\cite{GHJ}, 2.1) we have this way the following example.

Let $k=\mathbb{C}$ be the field of complex numbers, 
$\lambda^{-1} = 4\cos^2\frac{\pi}{n+3}\ (n\geq 2)$,
and $e_1, e_2,\dots$ be a sequence of idempotents satisfying,
for all $i$ and $j$,  the braid-like relations
\begin{eqnarray*}
e_i e_{i\pm 1} e_i &=& \lambda e_i, \\
e_i e_j &=& e_j e_i, \quad \mbox{if } |i-j| \geq 2.
\end{eqnarray*}
Let $A_{k,l}$ be the algebra generated by
$1, e_k, e_{k+1},\dots e_l$ ($k\leq l$), $\sigma$ be the algebra 
anti-endomorphism of $H= A_{1,2n-1}$ determined by $\sigma(e_i) = e_{2n-i}$
and $P_{k}\in A_{2n-k, 2n-1} \otimes A_{1,k}$ be the image of the 
separability idempotent of $A_{1,k}$ under $(\sigma\otimes \id)$.

Finally,  we denote by $\tau$ the non-degenerate Markov trace (\cite{GHJ}, 2.1)
on $H$ and by $w$ the index of the restriction of $\tau$ on
$A_{n+1, 2n-1}$ \cite{W}, i.e., the unique central element in $A_{n+1, 2n-1}$
such that $\tau(w\,\cdot )$ is equal to the trace of the left regular 
representation of  $A_{n+1, 2n-1}$.

Then the  following operations give a quantum groupoid structure on $H$ :
\begin{eqnarray*}
\Delta(yz) &=& (z\otimes y)P_{n-1}, \qquad y\in A_{n+1, 2n-1},\quad
               z\in A_{1,n-1}\\
\Delta(e_n) &=& (1\otimes w) P_{n} (1\otimes w^{-1}), \\
S(h)      &=& w^{-1}\sigma(h)w, \\
\eps(h)   &=& \lambda^{-n}\tau(hfw), \quad h\in A,  
\end{eqnarray*}
where in the last line
$$
f= \lambda^{n(n-1)/2}(e_ne_{n-1}\cdots e_1)(e_{n+1}e_n\cdots e_2)\cdots
   (e_{2n-1}e_{2n-2}\cdots e_n)
$$
is the Jones projection corresponding to the $n$-step basic construction.

The source and target counital subalgebras of $H= A_{1,2n-1}$  are 
$H_s=A_{n+1, 2n-1}$ and $H_t=A_{1,n-1}$ respectively. 
The example corresponding to $n=2$ is a $C^*$-quantum 
groupoid of dimension $13$ with the antipode having an infinite order
(it was studied in detail in \cite{NV2}, 7.3).

We refer the reader to Sections $8$ and $9$ of this survey 
and to the Appendix of \cite{NV3} for the explanation of how
quantum groupoids can be constructed from subfactors.

\subsection{Other versions of a quantum groupoid}
Here we briefly discuss several notions of quantum groupoids
that appeared in the literature and relations between them.
All these objects  generalize both usual groupoid algebras and
their duals and Hopf algebras. We apologize for possible
non-intentional omissions in the list below.

{\em Face algebras} of Hayashi \cite{H1} were defined as Hopf-like
objects containing an abelian subalgebra generated by ``bases''.
Non-trivial examples of such objects and applications to
monoidal categories and II${}_1$ subfactors were considered
in \cite{H2} and \cite{H3}. It was shown in (\cite{N}, 5.2) that
face algebras are precisely quantum groupoids whose counital
subalgebras are abelian.
 
{\em Generalized Kac algebras} of Yamanouchi \cite{Y} were used to
characterize $C^*$-algebras arising from finite groupoids; they are
exactly $C^*$-quantum groupoids with $S^2=\id$, see (\cite{NV1}, 2.5)
and (\cite{N}, 8.7).

The idea of a quantum transformation groupoid (see 2.6) has been explored 
in a different form in \cite{Mal1} (resp., \cite{V}), where it was shown 
that an action of a Hopf algebra on a commutative
(resp., noncommutative)
algebra gives rise to a specific quantum groupoid structure on the tensor 
product of these algebras. This construction served as a strong motivation
for \cite{Mal2}, where a quite general approach to quantum groupoids has 
been developed.

The definition of {\em quantum groupoids} in \cite{Lu} and \cite{Xu}
is more general than the one we use. In their approach
the existence of the bases is a part of the axioms, and they do not
have to be finite-dimensional. These objects give rise to Lie bialgebroids
as classical limits. It was shown in \cite{EN} that every quantum groupoid
in our sense is a quantum groupoid in the sense of \cite{Lu} and \cite{Xu},
but not vice versa.

A suitable functional analytic framework for studying quantum groupoids
(not necessarily finite) in the spirit of the S.~Baaj-G.~Skandalis 
multiplicative unitaries \cite{BS} was developed in \cite{EVal}, \cite{E2} 
in connection with depth 2 inclusions of von Neumann algebras. 
A closely related notion of a {\em Hopf bimodule} was introduced 
in \cite{V1} and then studied extensively in \cite{EVal}.
In the finite-dimensional case the equivalence of 
these notions to that of a $C^*$-quantum groupoid was shown in 
\cite{V2} and \cite{BSz3} and, respectively, in \cite{V2}, \cite{NV1}.

\end{section}

\begin{section}
{Integrals and semisimplicity}

\subsection{Integrals in quantum groupoids}
\begin{definition}[\cite{BNSz},  3.1]
\label{integral}
A left (right) {\em integral} in $H$ is an element
$l\in H$ ($r\in H$) such that 
\begin{equation}
hl =\eps_t(h)l, \qquad (rh = r\eps_s(h)) \qquad \mbox{ for all } h\in H. 
\end{equation}
 
\end{definition}
These notions clearly generalize the corresponding notions for Hopf 
algebras (\cite{M}, 2.1.1). We denote $\int_H^l$ (respectively, 
$\int_H^r$) the space of left (right) integrals in $H$ and by 
$\int_H = \int_H^l \cap \int_H^r$ the space of two-sided integrals.  

An integral in $H$ (left or right) is called {\em non-degenerate} if
it defines a non-degenerate functional on $\widehat H$. A left integral $l$
is called {\em normalized} if $\eps_t(l)=1$. Similarly, $r\in \int_H^r$ 
is normalized if $\eps_s(r)=1$. 

A dual notion to that of left (right) integral is the left 
(right) invariant measure. Namely, a functional $\phi \in \widehat H$ is said 
to be a left (right) {\em invariant measure} on $H$ if 
\begin{equation}
(\id\otimes \phi)\Delta = (\eps_t\otimes \phi)\Delta,\qquad  (\mbox{resp.}, 
(\phi\otimes \id)\Delta = (\phi\otimes \eps_s)\Delta).
\end{equation}
A left (right) invariant measure is said to be normalized if 
$(\id\otimes\phi)\Delta(1)=1$ (resp., $(\phi\otimes \id)\Delta(1)=1$).

\begin{example}
\label{examples of integrals}
(i) Let $G^0$ be the set of units of a finite groupoid $G$,
then the elements $l_e = \sum_{gg^{-1}=e}\,g\,(e\in G^0)$
span $\int_{kG}^l$ and elements $r_e = \sum_{g^{-1}g=e}\,g\,(e\in G^0)$
span $\int_{kG}^r$.
\newline
(ii) If $H = (kG)^*$ then $\int_H^l = \int_H^r =\span\{p_e,\,e\in G^0\}$.
\end{example}

The next proposition gives a description of the set of left integrals.
\begin{proposition}[\cite{BNSz}, 3.2]
\label{integral properties}
The following conditions for $l\in H$ are equivalent:
\begin{enumerate}
\item[(i)] $l\in \int_H^l$,
\item[(ii)] $(1\otimes h)\Delta(l) =(S(h)\otimes 1)\Delta(l)$ for all $h\in H$,
\item[(iii)] $(\id\otimes l)\Delta(\widehat H) =\widehat H_t$,
\item[(iv)] $(\Ker\,\eps_t)l =0$,
\item[(v)] $S(l)\in \int_H^r$. 
\end{enumerate}
\end{proposition}
\begin{proof}
The proof is a straightforward application of 
Definitions~\ref{finite quantum groupoid}, \ref{integral},
Propositions~\ref{counital properties} and \ref{S and counital maps},
and is left as an exercise for the reader.
\end{proof}

%
 
\subsection{Hopf modules}
Since a quantum groupoid $H$ is both algebra and coalgebra, one 
can consider modules and comodules over $H$. As in the theory of 
Hopf algebras, an $H$-Hopf module is an $H$-module which is also 
an $H$-comodule such that these two structures are compatible 
(the action ``commutes'' with coaction) :

\begin{definition}
\label{Hopf module}
A right $H$-Hopf module is such a $k$-vector space $M$ that
\begin{enumerate}
\item[(i)]
$M$ is a right $H$-module via $m\otimes h \mapsto m\cdot h$,
\item[(ii)]
$M$ is a right  $H$-comodule via $m\mapsto \rho(m) = m\0 \otimes m\I$,
\item[(iii)]
$(m\cdot h)\0 \otimes (m\cdot h)\I = m\0 \cdot h\1 \otimes m\I h\2$
\end{enumerate}
for all $m\in M$ and $h\in H$.
\end{definition}

\begin{remark}
Condition (iii) above means that $\rho$ is a right $H$-module map,
where $(M\otimes H)\Delta(1)$ is a right $H$-module via
$(m\otimes h)\Delta(1)\cdot h = (m\cdot h\1)\otimes (h\cdot h\2)$.
\end{remark} 

\begin{example}
$H$ itself is an $H$-Hopf module via $\rho =\Delta$.
\end{example}

\begin{example}
\label{dual Hopf module}
The dual vector space $\widehat H$
becomes a right $H$-Hopf module :
$$
\phi\cdot h = S(h)\act \phi \qquad 
\phi\0 \la \psi,\, \phi\I \ra =\psi\phi,
$$
for all $\phi,\psi\in \widehat H,\, h\in H$. Indeed, we need to check that
$$
(\phi\cdot h)\0 \otimes (\phi\cdot h)\I =
\phi\0\cdot h\1 \otimes \phi\1 h\2, \qquad \phi\in \widehat H, h\in H.
$$ 
Evaluating both sides against $\psi\in \widehat H$ in the first factor
we have :
\begin{eqnarray*}
(\phi\0\cdot h\1) \la \psi,\, \phi\1h\1 \ra
&=& (\psi\1 \phi)\cdot h\1 \la \psi\2,\, h\2 \ra \\
&=& \psi\1 \phi\1 \la S(\psi\2\phi\2),\, h\1 \ra \la \psi\3,\, h\2 \ra \\
&=& \psi\1 \phi\1 \la S^{-1}\eps_s(\psi\2)\phi\2,\,S(h) \ra \\
&=& \psi\phi\1 \la \phi\2,\, S(h) \ra \\ 
&=& \psi (S(h)\act \phi) = \psi(\phi\cdot h) \\
&=& (\phi\cdot h)\0 \la \psi,\, (\phi\cdot h)\I \ra.
\end{eqnarray*}
where we used Proposition~\ref{separability} and that
$S^{-1}\eps_s(\psi) \in \widehat H_t$ for all $\psi\in \widehat H$. 
\end{example}

The fundamental theorem for Hopf modules over Hopf algebras
([M], 1.9.4) generalizes to quantum groupoids as follows :

\begin{theorem}[\cite{BNSz}, 3.9]
\label{fundamental theorem}
Let $M$ be a right $H$-Hopf module
and 
\begin{equation}
N= Coinv M = \{m\in M \mid m\0\otimes m\I = m1\1 \otimes 1\2 \}
\end{equation}
be the set of coinvariants. The $H_t$-module tensor product $N\otimes_{H_t} H$
(where $N$ is a right $H_t$-submodule) is a right $H$-Hopf module via
\begin{equation}
(n\otimes h)\cdot g =n\otimes hg,
\qquad
(n\otimes h)\0 \otimes (n\otimes h)\I = (n\otimes h\1)\otimes h\2,
\end{equation}
for all $h,g\in H$ and $n\in N$. Then the map
\begin{equation}
\alpha : N\otimes_{H_t} H \to M : n\otimes h \mapsto n\cdot h
\end{equation}
is an isomorphism of right Hopf modules.
\end{theorem}
\begin{proof}
Proposition~\ref{counital properties}(iii) implies that 
the $H$-Hopf module structure on $M$ is well defined, and it is 
easy to check that $\alpha$ is a well defined homomorphism of 
$H$-Hopf modules. We will show that
\begin{equation}
\beta : M\to N\otimes_{H_t} H : m \mapsto (m\0 \cdot S(m\I))\otimes m\II
\end{equation}
is an inverse of $\alpha$. First, we observe that 
$m\0 \cdot S(m\I) \in N$ for all $m\in M$, since
\begin{eqnarray*}
\rho(m\0 \cdot S(m\I))
&=& (m\0 \cdot S(m\III)) \otimes m\I S(m\II) \\
&=& (m\0 \cdot S(1\2m\I) \otimes S(1\1) \\
&=& (m\0 \cdot S(m\I)1\1) \otimes 1\2,
\end{eqnarray*}
so $\beta$ maps to $N\otimes_{H_t} H$. Next, we check that it is both
module and comodule map :
\begin{eqnarray*}
\beta(m\cdot h) 
&=& m\0 \cdot h\1 S(m\I h\2) \otimes m\II h\3 \\
&=& m\0 \cdot \eps_t(h\1) S(m\I) \otimes m\II h\2 \\
&=& m\0 \cdot S(m\I 1\1) \otimes m\II  1\2h =\beta(m) \cdot h, \\
\beta(m\0) \otimes  m\I
&=& m\0 \cdot S(m\I) \otimes m\II \otimes m\III \\
&=& \beta(m)\0 \otimes \beta(m)\I,  
\end{eqnarray*}
Finally, we verify that $\alpha\circ \beta =\id$
and $\beta\circ\alpha = \id$ :
\begin{eqnarray*}
\alpha\circ \beta(m)
&=& m\0 \cdot S(m\I)m\II = m\0 \cdot \eps_s(m\I) \\
&=& m\0 \cdot 1\1\eps(m\I 1\2) = m\0 \eps(m\I) = m, \\
\beta\circ\alpha(n\otimes h)
&=& \beta(n\cdot h) = \beta(n) \cdot h \\
&=& n \cdot 1\1S(1\2) \otimes 1\3h = n \otimes h, 
\end{eqnarray*}
which completes the proof.
\end{proof}

\begin{remark}
\label{integrals are invariants}
$\int_{\widehat H}^l = \Coinv \widehat H$ where $\widehat H$ is 
an $H$-Hopf module as in Example~\ref{dual Hopf module}.
Indeed, the condition
$$
\lambda\I \otimes \lambda\II = \lambda\cdot 1\1 \otimes 1\2
$$
is equivalent to $\phi\lambda\ =(S(1\1)\act \lambda)\la \phi,\, 1\2\ra
= \eps_t(\phi)\lambda$.
\end{remark}

\begin{corollary}
\label{FT for integrals}
$\widehat H \cong \int_{\widehat H}^l \otimes_{H_t} H$ as right $H$-Hopf 
modules. In particular,  $\int_{\widehat H}^l$ is a non-zero subspace 
of $\widehat H$.
\end{corollary}
\begin{proof}
Take $M=\widehat H$ in Theorem \ref{fundamental theorem} and
use Remark~\ref{integrals are invariants}.
\end{proof}

\subsection{Maschke's Theorem}
The existence of left integrals in quantum groupoids
leads to the following generalization of Maschke's Theorem,  
well-known for Hopf algebras ([M], 2.2.1).

\begin{theorem}[\cite{BNSz}, 3.13]
\label{maschke}
Let $H$ be a finite quantum groupoid, then the following conditions
are equivalent :
\begin{enumerate}
\item[(i)] $H$ is semisimple,
\item[(ii)] There exists a normalized left integral $l$ in $H$, 
\item[(iii)] $H$ is separable. 
\end{enumerate}
\end{theorem}
\begin{proof}
(i) $\Rightarrow$ (ii) : Suppose that $H$ is semisimple, 
then since $\Ker\eps_t $ is a left ideal in $H$ we have 
$\Ker\eps_t =Hp$ for some idempotent $p$. Therefore,
$\Ker\eps_t (1-p) =0$ and $l=1-p$ is a left integral
by Lemma~\ref{integral properties}(iv). It is normalized since
$\eps_t(l) = 1-\eps_t(p) =1$. (ii) $\Rightarrow$ (iii) :
if $l$ is normalized then $l\1 \otimes S(l\2)$ is a separability
element of $H$ by Lemma~\ref{integral properties}(ii). 
(iii) $\Rightarrow$ (i) : this is a standard result \cite{P}.
\end{proof}

\begin{corollary}
\label{kG is semisimple}
Let $G$ be a finite groupoid and for every $e\in G^0$,
where $G^0$ denotes the unit space of $G$, let
$|e|=\#\{g\in G\mid gg^{-1}=e\}$. Then $kG$ is semisimple
iff $|e|\neq 0$ in $k$ for all $e\in G^0$.
\end{corollary}
\begin{proof}
In the notation of Example~\ref{examples of integrals}(i)
the element $l=\sum_{e\in G^0}\, \frac{1}{|e|}l_e$ is a normalized
integral in $kG$.
\end{proof}

Given a left integral $l$, one can show (\cite{BNSz}, 3.18) that 
if there exists $\lambda\in \widehat H$ such that $\lambda\act l=1$, then it 
is unique, it is a left integral in $\widehat H$ and $l\act\lambda=\widehat 1$. 
Such a pair $(l,\lambda)$ is called a {\em dual pair} of left integrals. 
One defines dual pairs of right integrals in a similar way.

\end{section}


\begin{section} 
{Actions and smash products}

\subsection{Module and comodule algebras}

\begin{definition}
\label{module algebra}
An algebra $A$ is a (left) {\em $H$-module algebra} if $A$ is a left
$H$-module via $h\otimes a \to h\cdot a$ and
\begin{enumerate}
\item[1)] $h\cdot ab = (h\1 \cdot a)(h\2 \cdot b)$,
\item[2)] $h\cdot 1 = \eps_t(h)\cdot 1$.
\end{enumerate} 
\end{definition}
If $A$ is an $H$-module algebra we will also say that $H$ acts on $A$.

\begin{definition}
\label{comodule algebra}
An algebra $A$ is a (right) {\em $H$-comodule algebra} if $A$ is a right
$H$-module via $\rho: a \to a^{(0)} \otimes a^{(1)}$ and
\begin{enumerate}
\item[1)] $\rho(ab) = a^{(0)} b^{(0)} \otimes a^{(1)} b^{(1)}$,
\item[2)] $\rho(1) = (\id\otimes\eps_t)\rho(1)$.
\end{enumerate} 
\end{definition}

It follows immediately that $A$ is a left $H$-module algebra
if and only if $A$ is a right $\widehat H$-comodule algebra.

\begin{example}
\label{examples of actions}
\begin{enumerate}
\item[(i)] The target counital subalgebra
$H_t$ is a trivial $H$-module algebra via
$h\cdot z = \eps_t(hz)$, $h\in H,\,z\in H_t$.
\item[(ii)] $H$ is an $\widehat H$-module algebra via
the dual action $\phi\act h = h\1 \la \phi,\,h\2\ra$, 
$\phi\in \widehat H,\,h\in H$.
\item[(iii)] Let $A = C_H(H_s) =\{ a\in H \mid ay=ya \quad\forall y\in H_s\}$,
be the centralizer of $H_s$ in $H$, then $A$ is an $H$-module algebra via
the adjoint action $ h\cdot a = h\1 a S(h\2)$.
\end{enumerate}
\end{example}

\subsection{Smash products}

Let $A$ be an $H$-module algebra, then a {\em smash product}
algebra $A\# H$ is defined on a $k$-vector space $A\otimes_{H_t} H$,
where $H$ is a left $H_t$-module via
multiplication and $A$ is a right $H_t$-module via
$$
a\cdot z = S^{-1}(z)\cdot a = a(z\cdot 1), \qquad a\in A, z\in H_t,
$$
as follows. Let $a\# h$ be the class of $a\otimes h$ in $A\otimes_{H_t} H$,
then the multiplication in $A\#H$ is given by the familiar formula
$$
(a\# h)(b\# g) = a(h\1 \# b) \# h\2 g,\qquad a,b,\in A,\, h,g\in H,
$$
and the unit of $A\# H$ is $1\# 1$.

\begin{example}
\label{trivial smash product}
$H$ is isomorphic to the trivial smash product algebra $H_t\# H$.
\end{example}

\subsection{Duality for actions}
An analogue of the Blattner-Montgomery duality theorem for
actions of quantum groupoids was proven in \cite{N2}.
Let $H$ be a finite quantum groupoid and
$A$ be a left $H$-module algebra. Then the smash product $A\#H$
is a left $\widehat H$-module algebra via
\begin{equation}
\label{dual action}
\phi\cdot (a\#h) = a\#(\phi\act h), \quad
\phi\in \widehat H,\,h\in H,\,a\in A.
\end{equation}
In the case when $H$ is a finite dimensional Hopf algebra,
there is an isomorphism
$(A\#H)\#\widehat H \cong M_n(A)$, where $n=\dim H$ and $M_n(A)$ is an algebra
of $n$-by-$n$ matrices over $A$ \cite{BM}.
This result extends to quantum groupoid
action in the form $(A\#H)\#\widehat H \cong \End(A\#H)_A$, 
where $A\#H$ is a right $A$-module
via multiplication (note that $A\#H$ is not necessarily a free $A$-module,
so that we have $\End(A\#H)_A\not\cong M_n(A)$ in general).

\begin{theorem}[\cite{N2}, 3.1, 3.2]
\label{BM-duality}
The map $\alpha : (A\#H)\#\widehat H \to \End(A\#H)_A$ defined by
$$
\alpha((x\#h)\#\phi)(y\#g) = (x\#h)(y\#(\phi\act g)) 
= x(h\1\cdot y) \#h\2 (\phi\act g) 
$$ 
for all $x,y\in A,\, h,g\in H,\,\phi\in \widehat H$ is an isomorphism of algebras.
\end{theorem}
\begin{proof}
A straightforward (but rather lengthly) computation shows  
that $\alpha$ is a well defined homomorphism, cf.\ (\cite{N2}, 3.1).
Let $\{ f_i\}$ be a basis of $H$ and $\{ \xi^i\}$ be the dual basis
of $\widehat H$, i.e., such that $\la f_i,\, \xi^j \ra =\delta_{ij}$ for
all $i,j$, then  the element $\sum_i\, f_i\otimes \xi^i \in H\otimes_k 
\widehat H$ does not depend on the choice of $\{ f_i\}$ and
the following map 
\begin{eqnarray*}
\beta 
&:& \End(A\# H)_A \to (A\# H)\# \widehat H \\
&:& T \mapsto \sum_i\, T(1\# {f_i}\2)(1\# S^{-1}({f_i}\1)) \# \xi^i.
\end{eqnarray*}
is the inverse of $\alpha$.
\end{proof}

\begin{corollary}
$H\# \widehat H \cong \End(H)_{H_t}$, therefore, $H\# \widehat H$ 
is a semisimple algebra.
\end{corollary}
\begin{proof} 
Theorem~\ref{BM-duality}
for $A=H_t$ shows that $H$ is a projective generating $H_t$-module
such that $\End(H)_{H_t} \cong H\# \widehat H$. Therefore, $H_t$ and
$H\# \widehat H$ are Morita equivalent. Since $H_t$ is semisimple
(as a separable algebra), $H\# \widehat H$ is semisimple.
\end{proof}

\end{section}


\begin{section}
{Representation category of a quantum groupoid}

Representation categories of quantum groupoids were studied in \cite{BSz2}
and \cite{NTV}.

\subsection{Definition of $\Rep(H)$}
For a quantum groupoid $H$
let $\Rep(H)$ be the category of  representations of $H$, whose objects are
$H$-modules of finite rank and whose morphisms are $H$-linear homomorphisms.
We show that, as in the case of Hopf algebras, $\Rep(H)$ has a natural 
structure of a monoidal category with duality.

For objects $V,W$ of $\Rep(H)$ set
\begin{equation}
V\otimes W = \{x \in V\otimes_k W \mid x= \Delta(1)\cdot x\},
\end{equation}
with the obvious action of $H$ via the comultiplication $\Delta$
(here $\otimes_k$ denotes the usual tensor product of vector spaces).

Since $\Delta(1)$ is an idempotent, 
$V\otimes W=\Delta(1)\cdot (V\otimes_k W)$.
The tensor product of morphisms is the restriction of usual
tensor product of
homomorphisms. The standard associativity isomorphisms
$\Phi_{U,V,W} : (U\otimes V)\otimes W \to U\otimes (V\otimes W)$
are functorial and satisfy the pentagon condition, since $\Delta$ is
coassociative. We will suppress these isomorphisms and write simply
$U\otimes V\otimes W$.

The target counital subalgebra $H_t\subset H$ has an $H$-module structure  
given by $h\cdot z = \eps_t(hz)$, where $h\in H,\, z\in H_t$.

\begin{lemma}
$H_t$ is the unit object of $\Rep(H)$.
\end{lemma}
\begin{proof}
Define a left unit homomorphism $l_V : H_t\otimes V \to V$ by
$$
l_V( 1\1\cdot z \otimes 1\2\cdot v) = z\cdot v, \qquad z\in H_t,\, v\in V.
$$
It is an invertible $H$-linear map  with the inverse 
$l_V^{-1}(v) = S(1\1) \otimes 1\2 \cdot v$. Moreover, the collection
$\{l_V\}_V$ gives a natural equivalence between the functor $H_t\otimes (\ )$
and the identity functor. 
Similarly, the right unit homomorphism
$r_V : V\otimes H_t \to V$ defined by
$$
r_V( 1\1\cdot v \otimes 1\2\cdot z) = S(z)\cdot v, \qquad z\in H_t,\, v\in V,
$$
has the inverse $r_V^{-1}(v) = 1\1 \cdot v \otimes 1\2$
and satisfies the necessary properties.
Finally, one can check the triangle axiom, i.e., that 
\begin{equation*}
(\id_V\otimes l_W) = (r_V\otimes id_W)
\end{equation*}
for all objects $V,W$ of $\Rep(H)$ and $v\in V,\,w\in W$ (\cite{NTV}, 4.1).
\end{proof}

Using the antipode $S$ of $H$, we can provide $\Rep(H)$ with a duality.
For any object $V$ of $\Rep(H)$ define the action of $H$ on $V^*=
\Hom_k(V,\,k)$ by $(h\cdot \phi)(v) = \phi(S(h)\cdot v)$, where $h\in H,
v\in V, \phi\in V^*$. For any morphism $f: V\to W$ let $f^* : W^*\to V^*$
be the morphism dual to $f$ (see \cite{T}, I.1.8).

For any $V$ in $\Rep(H)$ define the duality homomorphisms
\begin{equation*}
d_V : V^*\otimes V \to H_t, \qquad b_V: H_t \to V \otimes V^*
\end{equation*}
as follows. For $\sum_j\, \phi^j \otimes v_j\in V^* \otimes V$ set
\begin{equation}
d_V(\sum_j\, \phi^j \otimes v_j)= \sum_j\, \phi^j(1\1\cdot v_j)1\2.
\end{equation}
Let $\{g_i\}_i$ and $\{\gamma^i\}_i$ be bases of $V$ and $V^*$ respectively,
dual to each other. The element $\sum_i\,g_i\otimes \gamma^i$ does not
depend on choice of these bases; moreover, for all $v\in V, \phi\in V^*$
one has $\phi = \sum_i\,\phi(g_i)\gamma^i$ and $v =\sum_i\,g_i\gamma^i(v)$.
Set
\begin{equation}
b_V(z) = z\cdot \sum_i\, g_i \otimes \gamma^i.
\end{equation}
 
\begin{proposition}
\label{monoidal with duality}
The category $\Rep(H)$ is a monoidal category with duality.
\end{proposition}
\begin{proof}
We know already that $\Rep(H)$ is monoidal.
One can check (\cite{NTV}, 4.2) that $d_V$ and $b_V$ are $H$-linear.
To show that they satisfy the identities
$$
(\id_V\otimes d_V)(b_V \otimes \id_V) = \id_V, \qquad
(d_V \otimes \id_{V^*})(\id_{V^*}\otimes b_V) = \id_{V^*},
$$ 
take $\sum_j\, \phi^j \otimes v_j\in V^*\otimes V, z\in H_t$.
Using the isomorphisms $l_V$ and $r_V$ identifying $H_t\otimes V$,
$V\otimes H_t$ and $V$, for all $v\in V$ and $\phi\in V^*$ we have:
\begin{eqnarray*}
(\id_V\otimes d_V)(b_V \otimes \id_V)(v)
&=& (\id_V\otimes d_V)(b_V(1\1\cdot 1)\otimes 1\2\cdot v) \\
&=& (\id_V\otimes d_V)(b_V(1\2)\otimes S^{-1}(1\1)\cdot v) \\
&=& \sum_i\,(\id_V\otimes d_V)
            (1\2\cdot g_i \otimes 1\3\cdot \gamma^i \otimes S^{-1}(1\1)\cdot v) \\
&=& \sum_i\,1\2\cdot g_i \otimes (1\3\cdot \gamma^i)(1\1'S^{-1}(1\1)\cdot v)1\2'\\
&=& 1\2 S(1\3) 1\1' S^{-1}(1\1)\cdot v \otimes 1\2' =v, \\
(d_V \otimes \id_{V^*})(\id_{V^*}\otimes b_V)(\phi)
&=& (d_V \otimes \id_{V^*})(1\1\cdot \phi \otimes b_V(1\2)) \\
&=& \sum_i\, (d_V \otimes \id_{V^*})
             (1\1\cdot\phi \otimes 1\2\cdot g_i \otimes 1\3\cdot \gamma^i) \\
&=& \sum_i\,(1\1\cdot\phi)(1\1'1\2\cdot g_i)1\2' \otimes 1\3\cdot \gamma^i\\
&=& 1\2' \otimes 1\3 1\1 S(1\1'1\2)\cdot \phi = \phi,
\end{eqnarray*}
which completes the proof.
\end{proof}
\subsection{Quasitriangular quantum groupoids}

\begin{definition}
\label{QT WHA}
A quasitriangular quantum groupoid is a pair 
($H,\, \R$) where $H$ is a quantum groupoid and
$\R\in \Delta^{op}(1)(H\otimes_k H)\Delta(1)$
satisfying the following conditions :
\begin{equation}
\Delta^{op}(h)\R = \R\Delta(h),
\end{equation}
for all $h\in H$,
where $\Delta^{op}$ denotes the comultiplication opposite to $\Delta$,
\begin{eqnarray}
(\id \otimes \Delta)\R &=& \R_{13}\R_{12}, \label{DR1}\\
(\Delta \otimes \id)\R &=& \R_{13}\R_{23}, \label{DR2}
\end{eqnarray}
where $\R_{12} = \R\otimes 1$, $\R_{23} = 1\otimes \R$, etc.\ as usual,
and such that there exists $\bR\in \Delta(1)(H\otimes_k H)\Delta^{op}(1)$
with
\begin{equation}
\R\bR = \Delta^{op}(1), \qquad \bR\R = \Delta(1).
\end{equation}
\end{definition}

Note that $\bR$ is uniquely determined by $\R$: if $\bR$ and $\bR'$ are two
elements of $\Delta(1)(H\otimes_k H)\Delta^{op}(1)$ satisfying the previous
equation, then 
$$
\bR=\bR\Delta^{op}(1)=\bR\R\bR'=\Delta(1)\bR'=\bR'.
$$
For any two objects $V$ and $W$ of $\Rep(H)$ define
$c_{V,W}: V\otimes W \to W\otimes V$ as the action of $R_{21}$ :
\begin{equation}
c_{V,W}(x) = R\II\cdot x\II \otimes R\I\cdot x\I,
\end{equation}
where  $x =x\I\otimes x\II\in V\otimes W$ and  $\R=\R\I\otimes \R\II$.

\begin{proposition}
\label{braiding}
The family of homomorphisms $\{c_{V,W}\}_{V,W}$ defines a braiding in
$\Rep(H)$. Conversely, if $\Rep(H)$ is braided, then there exists $\R$,
satisfying the properties of Definition~\ref{QT WHA} and inducing
the given braiding.
\end{proposition}
\begin{proof}   
Note that $c_{V,W}$ is well-defined, since $\R_{21}= \Delta(1)\R_{21}$.
To prove the $H$-linearity of $c_{V,W}$ we observe that

\begin{eqnarray*}
c_{V,W}(h\cdot x)
&=& \R\II h\2\cdot x\II \otimes \R\I h\1\cdot x\I \\
&=& h\1 \R\II \cdot x\II \otimes h\2 \R\I \cdot x\I =h\cdot (c_{V,W}(x)).
\end{eqnarray*}
The inverse of $c_{V,W}$ is given by
$$
c_{V,W}^{-1}(y) = \bR^{(1)}\cdot y\II \otimes \bR^{(2)}\cdot y\I,
\qquad \mbox{ where } y =y\I\otimes y\II\in W\otimes V,
$$
therefore $c_{V,W}$ is is an isomorphism.
Finally, we check the braiding identities.
Let $ x =x\I\otimes x\II \otimes x\III \in U\otimes V\otimes W$, then
\begin{eqnarray*}
\lefteqn{(\id_V \otimes c_{U,W})(c_{U,V}\otimes \id_W)(x) = } \\
&=& (\id_V \otimes c_{U,W})
    (\R\II\cdot x\II \otimes \R\I\cdot x\I \otimes x\III \\
&=& \R\II\cdot x\II \otimes {\R'}\II\cdot x\III \otimes {\R'}\I \R\I\cdot x\I \\
&=& \R\II\1\cdot x\II \otimes \R\II\2\cdot x\III \otimes \R\I\cdot x\I
    = c_{U, V\otimes W}(x).
\end{eqnarray*}
Similarly, we have
$(c_{U,W}\otimes\id_V)(\id_U\otimes c_{V,W}) = c_{U\otimes V,W}$.

The third equality of this computation shows that
the relations of Definition~\ref{QT WHA}
are equivalent to the braiding identities.
\end{proof}

\begin{lemma}
\label{QYBE}
Let $(H,\R)$ be a quasitriangular quantum groupoid.
Then $\R$ satisfies the quantum Yang-Baxter equation :
$$
\R_{12}\R_{13}\R_{23} = \R_{23}\R_{13}\R_{12}.
$$
\end{lemma}  
\begin{proof}
It follows from the first two relations of Definition~\ref{QT WHA}, that
\begin{eqnarray*}
\R_{12}\R_{13}\R_{23}=(\id\otimes \Delta^{op})(\R) \R_{23}=
\R_{23} (\id\otimes \Delta)(\R)=\R_{23}\R_{13}\R_{12}.
\end{eqnarray*}
\end{proof}

\begin{remark}
Let us define two linear maps $\R_1, \R_2 : \widehat H \to H$ by
$$
\R_1(\phi) =(\id \otimes \phi)(\R), \quad
\R_2(\phi) =(\phi\otimes \id)(\R), \quad\mbox{for all } \phi\in \widehat H.
$$
Then condition (\ref{DR1}) of
Definition~\ref{QT WHA} is equivalent to $\R_1$
being a coalgebra homomorphism and algebra anti-homomorphism and
condition (\ref{DR2}) is equivalent
to $\R_2$ being an algebra homomorphism and coalgebra anti-homomorphism. 
In other words, $\R_1:\widehat H \to H^{\op}$ and 
$\R_2:\widehat H \to H^{\cop}$ are homomorphisms of quantum groupoids.
\end{remark}

\begin{proposition}
\label{properties of R}
For any quasitriangular quantum groupoid $(H,\R)$, we have: 
\begin{eqnarray*}
(\eps_s\otimes \id)(\R) &=& \Delta(1), \qquad
(\id\otimes \eps_s)(\R) = (S\otimes \id)\Delta^{op}(1), \\
(\eps_t\otimes \id)(\R) &=& \Delta^{op}(1),\qquad
(\id\otimes \eps_t)(\R) = (S\otimes \id)\Delta(1), \\
(S\otimes \id)(\R) &=& (\id\otimes S^{-1})(\R) = \bR,
\qquad (S\otimes S)(\R)= \R. 
\end{eqnarray*}
\end{proposition}
\begin{proof}
The proof is essentially the same as (\cite{Ma}, 2.1.5).
\end{proof}

\begin{proposition}
\label{elements u and v}  
Let $(H,\R)$ be a quasitriangular quantum groupoid. Then 
$$
S^2(h) = uhu^{-1}
$$ 
for all $h\in H$, where $u=S(\R\II)\R\I$ is
an invertible element of $H$ such that
$$
u^{-1} = \R\II S^2(\R\I), \quad \Delta(u) = \bR \bR_{21} (u\otimes u).
$$
Likewise, $v = S(u) = \R\I S(\R\II)$ obeys $S^{-2}(h) = vhv^{-1}$ and
$$
v^{-1} = S^2(\R\I)\R\II, \quad \Delta(v) = \bR \bR_{21} (v\otimes v).
$$
\end{proposition}
\begin{proof}
Note that $S(\R\II)y\R\I = S(y)u$ for all $y\in H_s$. Hence, we have
\begin{eqnarray*}
S(h\2)uh\1
&=& S(h\2) S(\R\II)\R\I h\1 = S(\R\II h\2)\R\I h\1 \\
&=& S(h\1 \R\II) h\2 \R\I = S(\R\II) \eps_s(h) \R\I \\
&=& S(\eps_s(h))u,
\end{eqnarray*}
for all $h\in H$. Therefore,
\begin{eqnarray*}
uh
&=& S(1\2)u1\2h = S(\eps_t(h\2)u h\1 \\
&=& S(h\2 S(h\3))u h\1 = S^2(h\3) S(h\2)u h\1 \\
&=&  S^2(h\2) S(\eps_s(h\1)) u = S(\eps_s(h\1)S(h\2))u = S^2(h)u.
\end{eqnarray*}
The remaining part of the proof follows the lines of
(\cite{Ma}, 2.1.8). The results for $v$ can be obtained
by applying the results for $u$ to the quasitriangular quantum groupoid
$(H^{op/cop}, \R)$.
\end{proof}

\begin{definition}
\label{Drinfeld element}
The element $u$ defined in Proposition~\ref{elements u and v} is
called {\em the Drinfeld element} of $H$.
\end{definition} 
  
\begin{corollary}
\label{uv is central}
The element $uv=vu$ is central and obeys
$$
\Delta(uv) =(\bR \bR_{21})^2(uv\otimes uv).
$$
The element $uv^{-1} =vu^{-1}$ is group-like and
implements $S^4$ by conjugation.
\end{corollary}   

\begin{proposition}
\label{range of F}
Given a quasitriangular quantum groupoid $(H, \R)$, consider a linear map
$F:\widehat H\to H$ given by

\begin{equation}
\label{map F}
F: \phi \mapsto (\phi\otimes \id)(\R_{21}\R),\qquad \phi\in \widehat H.
\end{equation}
Then the range of $F$ belongs to $C_H(H_s)$, the centralizer of $H_s$.
\end{proposition}
\begin{proof}
Take $y\in H_s$. Then we have :
$$
\phi(\R\II {\R'}\I) \R\I {\R'}\II y =
\phi(\R\II y {\R'}\I) \R\I {\R'}\II = \phi(\R\II {\R'}\I) y \R\I {\R'}\II ,
$$
therefore $F(\phi)\in C_H(H_s)$, as required.
\end{proof}

\begin{definition}[cf. \cite{Ma}, 2.1.12]
\label{factorizability}
A quasitriangular quantum groupoid is {\em factorizable} if  
the above map $F: \widehat H \to C_H(H_s)$ is surjective.
\end{definition}

The factorizability of $H$ means that $\R$
is as non-trivial as possible, in contrast to {\em triangular} 
quantum groupoids, for which $\R_{21} = \bR $ and the range of $F$ is
equal to $H_t$.


\subsection{The Drinfeld double}

To define the {\em Drinfeld double} $D(H)$ of a 
quantum groupoid $H$, consider on the vector space 
$\widehat H^{op}\otimes_k H$ an associative multiplication
\begin{equation}
(\phi\otimes h)(\psi\otimes g) =
\psi\2\phi \otimes h\2 g \la S(h\1),\, \psi\1\ra \la h\3,\, \psi\3\ra,
\end{equation}
where $\phi, \psi \in \widehat H^{\op}$ and $h,g \in H$.
Then one can verify that the linear span $J$ of the elements
\begin{eqnarray}
\label{amalgamation}
\phi\otimes zh &-& (\eps\actl z)\phi \otimes h, \quad z\in H_t,\\
\phi\otimes yh &-& (y\actr\eps)\phi \otimes h, \quad y\in H_s,
\end{eqnarray}
is a two-sided ideal in $\widehat H^{op}\otimes_k H$. Let $D(H)$ be the
factor-algebra $(\widehat H^{op}\otimes_k H)/J$ and let $[\phi\otimes h]$ 
denote the class of $\phi\otimes h$ in $D(H)$.

\begin{proposition}
\label{the double}
$D(H)$ is a quantum groupoid with the unit $[\eps\otimes 1]$,
and comultiplication, counit, and antipode given by
\begin{eqnarray}
\Delta([\phi\otimes h])
&=& [\phi\1 \otimes h\1] \otimes [\phi\2 \otimes h\2],\\
\eps([\phi\otimes h])
&=& \la \eps_t(a),\,\phi \ra,\\
S([\phi\otimes a])
&=& [ S^{-1}(\phi\2) \otimes S(h\2)] \la h\1,\,\phi\1 \ra 
   \la S(h\3),\,\phi\3\ra.
\end{eqnarray}
\end{proposition}
\begin{proof}
The proof is a straightforward verification that all the structure
maps are well-defined and satisfy the axioms of a quantum groupoid,
which is carried out in full detail in (\cite{NTV}, 6.1).
\end{proof}

\begin{proposition}
\label{QT of D(H)}
The Drinfeld double $D(H)$ has a canonical quasitriangular structure
given by
\begin{equation}
\label{R-matrix}
\R = \sum_i [\xi^i \otimes 1] \otimes [\eps \otimes f_i],
\qquad
\bR  = \sum_j\, [S^{-1}(\xi_j) \otimes 1] \otimes[\eps \otimes f_j]
\end{equation}
where $\{ f_i\}$ and $\{ \xi^i\}$ are bases in $H$ and $\widehat H$
such that $\la f_i,\, \xi^j\ra =\delta_{ij}$.
\end{proposition}
\begin{proof}
The identities $(\id \otimes \Delta)\R = \R_{13}\R_{12}$ and
$(\Delta \otimes \id)\R = \R_{13}\R_{23}$
can be written as (identifying $[(\widehat H)^{op}\otimes 1]$ with $(\widehat H)^{op}$ 
and $[\eps \otimes H]$ with $H$) :
\begin{eqnarray*}
\sum_i\, \xi\1^i \otimes \xi\1^i \otimes f_i &=&
\sum_{ij}\, \xi^i \otimes \xi^j \otimes f_i f_j, \\
\sum_i\, \xi^i \otimes {f_i}\1  \otimes {f_i}\2 &=&
\sum_{ij}\, \xi^j \xi^i \otimes f_j \otimes f_i.
\end{eqnarray*}
The above equalities can be verified, e.g., by evaluating both sides
against an element $a\in H$ in the third factor (resp.,\ against
$\phi\in(\widehat H)^{op}$ in the first factor), see (\cite{Ma}, 7.1.1).
It is also straightforward to check that $\R$ is an intertwiner between
$\Delta$ and $\Delta^{op}$ and check that $\bR \R = \Delta(1)$ and
$\R\bR = \Delta^{\op}(1)$, see (\cite{NTV}, 6.2).
\end{proof}

\begin{remark}
\label{dual double}
The dual quantum groupoid $\widehat{D(H)}$ consists of all 
$\sum_k\, h_k\otimes \phi_k$ in $H\otimes_k {\widehat H}^{op}$ such that
$$
\sum_k\, (h_k\otimes \phi_k)|_J = 0.
$$
The structure operations of $\widehat{D(H)}$ are obtained by dualizing those
of $D(H)$:
\begin{eqnarray*}  
(  \sum_k\, h_k\otimes \phi_k ) (  \sum_l\, g_l\otimes \psi_l )
&=& \sum_{kl}\, h_kg_l\otimes \phi_k \psi_l, \\
1_{\widehat{D(H)}}
&=& 1\2 \otimes (\eps \actl 1\1), \\
\Delta(\sum_k\, h_k\otimes \phi_k )
&=& \sum_{ijk}\, ({h_k}\2 \otimes \xi^i {\phi_k}\1 \xi^j) \otimes
    (S(f_i) {h_k}\1 f_j \otimes {\phi_k}\2), \\
\eps(  \sum_k\, h_k\otimes \phi_k )
&=& \sum_k\,\eps(h_k)\widehat\eps(\phi_k),\\
S(  \sum_k\, h_k\otimes \phi_k )
&=& \sum_{ijk}\, f_i S^{-1}(h_k)S(f_j) \otimes \xi^i S(\phi)\xi^j,
\end{eqnarray*}
for all $\sum_k\, h_k\otimes \phi_k,\, \sum_l\, g_l\otimes \psi_l \in
\widehat{D(H)}$, where $\{ f_i\}$ and $\{\xi^j\}$ are dual bases.
\end{remark}

\begin{corollary}[\cite{NTV}, 6.4]
\label{D(H) is factorizable}
The Drinfeld double $D(H)$ is factorizable in the sense of
Definition~\ref{factorizability}.
\end{corollary}
\begin{proof}
One can use the explicit form of the $R$-matrix (\ref{R-matrix}) of $D(H)$
and the description of the dual $\widehat{D(H)}$ to check that  in this case
the map $F$ from (\ref{map F}) is surjective.
\end{proof}

\subsection{Ribbon quantum groupoids}

\begin{definition}
\label{Ribbon WHA}
A ribbon quantum groupoid is a quasitriangular quantum groupoid 
with an invertible central element $\nu\in H$ such that
\begin{equation}
\Delta(\nu) = \R_{21}\R(\nu\otimes \nu) \quad\mbox{and}\quad S(\nu)=\nu.
\end{equation}
The element $\nu$ is called a {\em ribbon element} of $H$.
\end{definition}

For an object $V$ of $\Rep(H)$ we define  the twist $\theta_V :V\to V$
to be the multiplication by $\nu$ :
\begin{equation}
\theta_V(v) =\nu\cdot v, \quad v\in V.
\end{equation}

\begin{proposition}
\label{twist}
Let $(H, \R,\nu)$ be a ribbon quantum groupoid. The family of
homomorphisms $\{\theta_V\}_V$ defines a twist in the braided monoidal
category $\Rep(H)$ compatible with duality.
Conversely, if $\theta_V(v) =\nu\cdot v$ is a twist in $\Rep(H)$,
then $\nu$  is a ribbon element of $H$.
\end{proposition}
\begin{proof}
Since $\nu$ is an invertible central element of $H$, the homomorphism
$\theta_V$ is an $H$-linear isomorphism. The twist identity follows
from the properties of $\nu$ :
$$
c_{W,V}c_{V,W}(\theta_V\otimes\theta_W)(x) =
\R_{21}\R(\nu\cdot x\I \otimes \nu\cdot x\II) = \Delta(\nu)\cdot x
= \theta_{V\otimes W}(x),
$$
for all $x =x\I\otimes x\II\in V\otimes W$.
Clearly, the identity $\R_{21}\R(\nu\otimes \nu) =\Delta(\nu)$ is equivalent
to the twist property. It remains to prove that
$$
(\theta_V \otimes \id_{V^*})b_V(z) = (\id_V \otimes \theta_{V^*})b_V(z),
$$
for all $z\in H_t$, i.e., that
$$
\sum_i\,\nu z\1\cdot \gamma^i \otimes z\2\cdot g_i =
\sum_i\,z\1\cdot \gamma^i \otimes \nu  z\2\cdot g_i,
$$
where $\sum_i\,\gamma^i \otimes g_i$ is the canonical element
in $V^*\otimes V$. Evaluating the first factors of the above equality
on an arbitrary $v\in V$, we get the equivalent condition :
$$
\sum_i\,(\nu z\1\cdot \gamma^i)(v) z\2\cdot g_i =
\sum_i\, (z\1\cdot \gamma^i)(v) \nu z\2\cdot g_i,
$$
which reduces to $z\2S(\nu z\1)\cdot v = S(z\1)\nu z\2\cdot v$.
The latter easily follows from the centrality of $\nu=S(\nu)$ and
properties of $H_t$.
\end{proof}

\begin{proposition}
\label{ribbon category}
The category $\Rep(H)$ is a ribbon category if and only if  
$H$ is a ribbon quantum groupoid.
\end{proposition}
\begin{proof}
Follows from   Propositions~\ref{monoidal with duality}, \ref{braiding},
and \ref{twist}.
\end{proof}
 
For any endomorphism $f$ of the object $V$, we define, following 
\cite{T}, I.1.5, its {\em quantum trace}
\begin{equation}
\label{q1-trace}
\tr_q(f)=d_V c_{V,V^*}(\theta_Vf \otimes \id_{V^*})b_V
\end{equation}
with values in $\End(H_t)$ and the {\em quantum dimension} of $V$
by $\dim_q(V)=\tr_q(id_V)$. 

\begin{corollary}[\cite{NTV}, 7.4]
\label{quantum trace}
Let $(H, \R,\nu)$ be a ribbon quantum groupoid, $f$ be an endomorphism
of an object $V$ in $\Rep(H)$. Then 
\begin{equation}
\tr_q(f)(z) = \Tr(S(1\1)u\nu f)z1\2,\qquad \dim_q(V)(z)=\Tr(S(1\1)u\nu)z1\2,
\end{equation}
where $\Tr$ is the usual trace of endomorphism, and $u$ is
the Drinfeld element.
\end{corollary}

\begin{corollary}
\label{quantum trace for connected WHA}
Let $k$ be algebraically closed.
If the trivial $H$-module $H_t$ is irreducible (which happens exactly
when $H_t \cap Z(H) = k$, i.e., when $H$ is connected (\cite{N1}, 3.11,
\cite{BNSz}, 2.4), then $\tr_q(f)$ and $dim_q(V)$ are scalars:
\begin{equation}
\label{q-trace}
\tr_q(f) = (\dim H_t)^{-1}\ \Tr(u\nu f),\qquad 
\dim_q(V) = (\dim H_t)^{-1} \Tr(u\nu).
\end{equation}
\end{corollary}
\begin{proof}
Any endomorphism of an irreducible module is the multiplication by
a scalar, therefore, we must have $\Tr(S(1\1)u\nu f)1\2 = \tr_q(f) 1$.
Applying the counit to both sides and using that $\eps(1) = \dim H_t$,
we get the result.
\end{proof}

\subsection{Towards modular categories}

In \cite{RT} a general method of constructing invariants of
$3$-manifolds from modular Hopf algebras was introduced. 
Later it became clear that the technique of Hopf algebras can
be replaced by a more general technique of
modular categories (see \cite{T}). In addition to quantum groups,
such categories also arise from skein categories of tangles
and, as it was observed by A.~Ocneanu, from certain bimodule
categories of type II${}_1$
subfactors.

The representation categories of {\em quantum groupoids} give
quite general construction of modular categories.
Recall some definitions needed here. 
Let ${\mathcal V}$ be a ribbon $Ab$-category over $k$, i.e.,
such that all $\Hom(V,W)$ are $k$-vector spaces (for all objects
$V,W\in {\mathcal V}$) and both operations $\circ$ and $\otimes$ are
$k$-bilinear.

An object $V\in {\mathcal V}$ is said to be {\it simple} if any 
endomorphism of $V$ is multiplication by an element of $k$.
We say that a family $\{V_i\}_{i\in I}$ of objects of ${\mathcal V}$ 
dominates  an object $V$ of ${\mathcal V}$ if there exists a finite set
$\{V_{i(r)}\}_r$ of objects of this family (possibly, with repetitions)
and a family of morphisms  $f_r:V_{i(r)} \to V,g_r:V\to V_{i(r)}$ such that 
$\id_V = \sum_r f_r g_r$.

A modular category (\cite{T}, II.1.4) is a pair consisting of
a ribbon $Ab$-category ${\mathcal V}$ and a finite family $\{V_i\}_{i\in I}$
of simple objects of ${\mathcal V}$ satisfying four axioms:
\begin{enumerate}
\item[(i)] There exists $0\in I$ such that $V_0$
is the unit object.
\item[(ii)] For any $i\in I$, there exists $i^*\in I$ such that $V_{i^*}$ is 
isomorphic to $V^*_i$.
\item[(iii)] All objects of ${\mathcal V}$ are dominated by the
family $\{V\}_{i\in I}$.
\item[(iv)] The square matrix $S=\{S_{ij}\}_{i,j\in I} =
\{ \tr_q(c_{V_i,V_j}\circ c_{V_j,V_i}) \}_{i,j\in I}$ is invertible over $k$
(here $\tr_q$ is the quantum trace in a ribbon category defined by 
(\ref{q1-trace})).
\end{enumerate}

If a quantum groupoid $H$ is connected and semisimple over an 
algebraically closed field, modularity of $\Rep(H)$ is equivalent to 
$\Rep(H)$ being ribbon and such that the matrix $S=\{S_{ij}\}_{i,j\in I} =
\{ \tr_q(c_{V_i,V_j}\circ c_{V_j,V_i}) \}_{i,j\in I}$, where $I$ is the set of 
all (equivalent classes of) irreducible representations, is invertible.
The following proposition extends a result known for Hopf algebras 
(\cite{EG}, 1.1).
\begin{proposition}
\label{factorizable implies modular}
If $H$ is a connected, ribbon, factorizable quantum groupoid 
over an algebraically closed field $k$,
possessing a normalized two-sided integral, then $\Rep(H)$ is a modular category.
\end{proposition}
\begin{proof}
Here we only need to prove the invertibility of the matrix formed by
\begin{eqnarray*}
S_{ij} 
&=& \tr_q(c_{V_i,V_j}\circ c_{V_j,V_i}) \\
&=& (\dim H_t)^{-1} \Tr( (u\nu)\circ c_{V_i,V_j}\circ c_{V_j,V_i}) \\
&=& (\dim H_t)^{-1}  (\chi_j \otimes \chi_i)( (u\nu\otimes u\nu)\R_{21}\R),
\end{eqnarray*}
where $V_i$ are as above, $I = \{1,\dots n\},\ \{ \chi_j\}$ is a basis 
in the space $C(H)$ of characters 
of $H$ (we used above the  formula (\ref{q-trace}) for the quantum trace).

It was shown in (\cite{NTV}, 5.12, 8.1) that the map $F : 
\phi\mapsto (\phi\otimes \id)(R_{21}R)$ is a linear isomorphism between
$C(H)\actl u\nu$ 
and the center $Z(H)$ (here $(\chi\actl u\nu)(a) =
\chi(u\nu a)~\forall a\in H,\ \chi\in C(H)$).
So, there exists an invertible matrix $T=(T_{ij})$ representing 
the map $F$ in the bases $\{ \chi_j\}$ of $C(H)$ and $\{e_i\}$ of $Z(H)$, 
i.e., such that
$F(\chi_j\actl u\nu) = \sum_i\, T_{ij}e_i$. Then
\begin{eqnarray*}
S_{ij} 
&=& (\dim H_t)^{-1} \chi_i(u\nu F(\chi_j\actl u\nu)) 
    = (\dim H_t)^{-1} \sum_k\, T_{kj} \chi_i(u\nu e_k) \\
&=& (\dim H_t)^{-1} (\dim V_i) \chi_i(u\nu) T_{ij}. 
\end{eqnarray*}
Therefore, $S = DT$, where $D = \mbox{diag}\{ (\dim H_t)^{-1} (\dim V_i) 
\chi_i(u\nu) \}$.
Theorem\ref{Haar theorem} below, shows that the existence of
a normalized two-sided integral in $H$ is equivalent to $H$
being semisimple and possessing an invertible element $g$ such
that $S^2(x)=gxg^{-1}$ for all $x\in H$ and $\chi(g^{-1})\neq O$
for all irreducible character $\chi$ of $H$. 
Then
$u^{-1}g$ is an invertible central element of $H$ and
$\chi_i(u^{-1})\neq 0$ for all $\chi_i$. By Corollary~\ref{uv is central} $uS(u)=c$
is invertible central, therefore 
$\chi_i(u)= \chi_i(c)\chi_i(S(u^{-1})))\neq 0$. 
Hence, $\chi_i(u\nu)\neq 0$ for all $i$ and $D$ is invertible.
\end{proof}
\end{section}


\begin{section}
{Twisting and dynamical quantum groups}

In this section we present a generalization of the Drinfeld 
twisting construction to quantum groupoids
developed in \cite{NV1}, \cite{Xu}, and \cite{EN}. We show that
{\em dynamical} twists of Hopf algebras give rise to quantum groupoids.
An important concrete example is given by dynamical quantum groups
at roots of $1$ \cite{EN}, which are dynamical deformations of the 
Drinfeld-Jimbo-Lusztig quantum groups $U_q(\g)$. 
The resulting quantum groupoids turn out to be selfdual, which is
a fundamentally new property, not satisfied by $U_q(\g)$. 

Most of the material of this section is taken from \cite{EN}.

\subsection{Twisting of quantum groupoids}

We describe the procedure of constructing new quantum groupoids
by twisting a comultiplication. Twisting of Hopf algebroids without
an antipode was developed in \cite{Xu} and  a special case of
twisting of $*$-quantum groupoids was considered in \cite{NV1}.

\begin{definition}
\label{a twist}
A {\em twist} for a quantum groupoid $H$ is a pair ($\Theta,\,\bTheta$), with
\begin{equation}
\label{Theta and bTheta}
\Theta\in \Delta(1)(H\otimes_k H), \quad
\bTheta\in (H\otimes_k H)\Delta(1), \quad
\mbox { and } \quad
\Theta\bTheta= \Delta(1)
\end{equation}
satisfying the following axioms :
\begin{eqnarray}
(\eps\otimes\id)\Theta = (\id\otimes \eps)\Theta 
&=& (\eps\otimes\id)\bTheta = (\id\otimes \eps)\bTheta = 1,
\label{eqn: twist eps} \\
(\Delta\otimes\id)(\Theta) (\Theta\otimes 1) &=&
(\id\otimes \Delta)(\Theta) (1\otimes\Theta),
\label{eqn: twist ++}\\
(\bTheta\otimes 1) (\Delta\otimes\id)(\bTheta) &=&
(1\otimes \bTheta) (\id\otimes \Delta)(\bTheta),
\label{eqn: twist --}\\ 
(\Delta\otimes\id)(\bTheta) (\id\otimes \Delta)(\Theta)  &=&
(\Theta\otimes 1) (1\otimes \bTheta),  
\label{eqn: twist +-}\\
(\id\otimes \Delta)(\bTheta) (\Delta\otimes\id)(\Theta) &=&
(1\otimes\Theta) (\bTheta\otimes 1).
\label{eqn: twist -+}
\end{eqnarray}  
\end{definition}

Note that for ordinary Hopf algebras the above notion coincides with the usual
notion of twist and each of the four conditions~(\ref{eqn: twist ++})
-- (\ref{eqn: twist -+}) implies the other
three. But since $\Theta$ and $\bTheta$ are, in general,
not invertible we need to impose all of them.

The next Proposition extends Drinfeld's  twisting construction
to the case of quantum groupoids.

\begin{proposition}
\label{twisting}
Let ($\Theta,\bTheta)$  be a twist for a quantum groupoid $H$.
Then there is a quantum groupoid $H_\Theta$ having
the same algebra structure and counit as $H$ with
a comultiplication and antipode given by
\begin{equation}
\Delta_\Theta(h) = \bTheta\Delta(h)\Theta, \qquad
S_\Theta(h) = v^{-1}S(h)v,
\label{eqn: theta}
\end{equation}
for all $h\in H_\Theta$, where $v= m(S\otimes\id)\Theta$
is invertible in $H_\Theta$.
\end{proposition}
\begin{proof}
The proof is a straightforward verification, see \cite{EN} for details.
\end{proof}

\begin{remark}
\label{rem: twist}
\begin{enumerate}
\item[(a)] One can check that the counital maps of the twisted quantum groupoid
$H_\Theta$ are given by
\begin{equation*}
(\eps_t)_\Theta(h) = \eps(\Theta\I h)\Theta\II, \qquad
(\eps_s)_\Theta(h) = \bTheta^{(1)}  \eps(h \bTheta^{(2)}).
\end{equation*}
\item[(b)] It is possible to  generalize the above twisting construction 
by weakening the counit condition~(\ref{eqn: twist eps}) and requiring only 
the existence of $u,w\in H$ such that 
\begin{equation}
\label{general eps twist}
\eps(\Theta\I w)\Theta\II = \Theta\I \eps(\Theta\II w) = 1,
\qquad
\eps(u \bTheta\I) \bTheta\II = \bTheta\I \eps(u \bTheta\II) =1.
\end{equation}
Then the counit of the twisted quantum groupoid $H_\Theta$ is 
\begin{equation}   
\eps_\Theta(h) = \eps(uhw), \qquad h\in H.
\end{equation}
\item[(c)] 
If $(\Theta,\bTheta)$ is a twist for $H$ and $x\in H$ is an invertible
element such that $\eps_t(x)=\eps_s(x)=1$ then $(\Theta^x,\bTheta^x)$, where
$$
\Theta^x =\Delta(x)^{-1}\Theta (x\otimes x) \quad \mbox{ and } \quad
\bTheta^x = (x^{-1}\otimes x^{-1})\bTheta \Delta(x),
$$
is also a twist for $H$. The twists $(\Theta,\bTheta)$ and
$(\Theta^x,\bTheta^x)$ are called {\em gauge equivalent} and
$x$ is called a {\em gauge transformation}.  Given such an $x$,
the map $h\mapsto x^{-1}hx$ is an isomorphism between quantum
groupoids $H_\Theta$ and $H_{\Theta^x}$.
\item[(d)] 
A twisting of a quasitriangular quantum groupoid is again
quasitriangular. Namely, if $(\Theta,\bTheta)$ is a twist
and $(\R,\bR)$ is a quasitriangular structure  for $H$ then
the quasitriangular structure for $H_\Theta$ is given by
$(\bTheta_{21}\R \Theta, \bTheta \bR \Theta_{21})$.
The proof of this fact is exactly the same as for Hopf algebras.
\end{enumerate}
\end{remark}

\begin{example}
\label{twisting of elementaries}
We will show that for every (non-commutative) separable algebra $B$
there is a family of quantum groupoid structures $H_q$ 
on the vector space $B^{op}\otimes B$ considered in (\cite{BSz2}, 5.2)
that can be understood in terms of 
twisting.  Let $H=H_{1_B}$ be  a quantum groupoid
with the following structure operations :
\begin{equation*}
\Delta(b\otimes c) = (b\otimes e\I) \otimes (e\II \otimes c), \qquad
\eps(b\otimes c) = \omega(bc), \qquad
S(b\otimes c) = c\otimes b,
\end{equation*}
for all $b,c\in B$, where $e =e\I \otimes e\II$ 
is the symmetric separability  idempotent of $B$ and $\omega$
is defined by the condition $(\omega\otimes\id)e =1$.
Note that $S^2=\id$ and $H$, as an algebra, is generated by 
its counital subalgebras, i.e., it is dual to an {\em elementary} 
quantum groupoid of (\cite{NV1}, 3.2).

Observe that the pair
$$
\Theta = (1\otimes e\I q^{-1}) \otimes (e\II \otimes 1)
\qquad
\bTheta = (1\otimes e\I) \otimes (e\II \otimes 1) = \Delta(1),
$$
where $q$ is an invertible element of $B$ with $e\I q e\II =1$,
satisfies (\ref{eqn: twist ++}) -- (\ref{eqn: twist -+})
and that conditions (\ref{general eps twist}) hold for $u=1$ and
$w=q$. According to Proposition~\ref{twisting} and Remark~\ref{rem: twist}(b),
these data define a twisting of $H$ such that the comultiplication and
antipode of the twisted quantum groupoid  are given by 
\begin{equation*}
\Delta(b\otimes c) = (b\otimes e\I q^{-1}) \otimes (e\II \otimes c), \quad
\eps(b\otimes c) = \omega(qbc),\quad
S(b\otimes c) = q^{-1}cq\otimes b,
\end{equation*}
for all $b,c\in B$.
The target and source counital subalgebras
of $H_q$ are $B^{op} \otimes 1$ and $1\otimes B$.     
The square of the antipode is
implemented by $g_q =q\otimes q$. Since $\Ad g_q$ is
an invariant of $H_q$, this example shows that there
can be uncountably many non-isomorphic semisimple 
quantum groupoids with the same underlying algebra (for
noncommutative $B$).

\end{example}

\subsection{Dynamical twists of Hopf algebras}

We describe a method of constructing twists of quantum groupoids,
which is a  finite-dimensional modification of the construction
proposed in \cite{Xu} and is dual to that of \cite{EV}, cf.\ \cite{EN}.

Dynamical twists  first appeared in the work of Babelon \cite{Ba},
see also \cite{BBB}.

Let $U$ be a Hopf algebra and $A=\Map(\mathbb{T},\, k)$
be a commutative and cocommutative Hopf algebra
of functions on a finite Abelian group $\mathbb{T}$
which is a Hopf subalgebra of $U$. Let $P_\mu,\, \mu\in \mathbb{T}$
be the minimal idempotents in $A$.

\begin{definition} 
\label{zero weight}
We say that an element $x$ in $U^{\otimes n},\, n\geq 1$
has {\em zero weight} if $x$ commutes with
$\Delta^n(a)$ for all $a\in A$, where $\Delta^n: A \to A^{\otimes n}$
is the iterated comultiplication.
\end{definition}

\begin{definition}
\label{dynamical twist}
An invertible, zero-weight  $U^{\otimes 2}$-valued function
$J(\lambda)$ on $\mathbb{T}$ is called a {\em dynamical twist }
for $U$ if it satisfies
the following functional equations :
\begin{eqnarray}
\label{dynamical equation}
(\Delta\otimes \id)J(\lambda) (J(\lambda+h^{(3)}) \otimes 1)
&=& (\id \otimes\Delta) J(\lambda) (1\otimes J(\lambda) ), \\
(\eps\otimes \id)J(\lambda) = (\id \otimes \eps)J(\lambda)  &=& 1.
\end{eqnarray}
Here an in what follows the notation $\lambda+ h^{(i)}$ means that   
the argument $\lambda$ is shifted by the weight of the $i$-th component,
e.g., $J(\lambda+h^{(3)}) =\sum_\mu\, J(\lambda+\mu) \otimes P_\mu
\in U^{\otimes 2} \otimes_k A$.
\end{definition}  

Note that for every fixed $\lambda\in \mathbb{T}$ the element
$J(\lambda)\in U\otimes U$ does not have to be a twist for $U$
in the sense of Drinfeld.
It turns out that an appropriate object for which $J$
defines a twisting is a certain quantum groupoid that
we describe next.

Observe that the simple algebra $\End_k(A)$ has a natural
structure of a cocommutative quantum groupoid given as follows
(cf.\ Section~\ref{groupoid algebras}) : 

Let $\{ E_{\lambda\mu} \}_{\lambda,\mu\in \mathbb{T}}$ be a basis
of $\End_k(A)$ such that
\begin{equation}
(E_{\lambda\mu}f)(\nu) = \delta_{\mu\nu}f(\lambda), \qquad
f\in A,\, \lambda,\mu,\nu\in \mathbb{T},
\end{equation}
then the comultiplication, counit, and antipode of
$\End_k(A)$ are given by
\begin{equation}
\Delta(E_{\lambda\mu}) = E_{\lambda\mu} \otimes E_{\lambda\mu},
\quad \eps(E_{\lambda\mu}) =1, \quad 
S(E_{\lambda\mu}) = E_{\mu\lambda}.
\end{equation}

Define the tensor product quantum groupoid $H = \End_k(A)\otimes_k U$
and observe that the elements
\begin{equation}
\Theta = \sum_{\lambda\mu}\, E_{\lambda\lambda+\mu} \otimes
E_{\lambda\lambda}P_\mu  \quad \mbox{ and } \quad
\bTheta = \sum_{\lambda\mu}\, E_{\lambda+\mu\lambda } \otimes
E_{\lambda\lambda}P_\mu
\label{eqn: Xu's theta}
\end{equation}
define a twist for $H$. Thus, according to Proposition~\ref{twisting},
$H_\Theta = (\End_k(A)\otimes_k U)_\Theta$ becomes
a quantum groupoid . It is non-commutative, non-cocommutative,
and not a Hopf algebra if $|\mathbb{T}|>1$.

It was shown in \cite{EN}, following \cite{Xu}, that $H_\Theta$ can be further
twisted by means of a dynamical twist $J(\lambda)$ on $U$. Namely,
if $J(\lambda)$ is a dynamical twists on $U$ embedded in $H\otimes_k H$ as 
\begin{equation}
J(\lambda) = \sum_{\lambda} E_{\lambda\lambda}
J\I(\lambda) \otimes E_{\lambda\lambda}J\II(\lambda),
\end{equation}
then  the pair ($F(\lambda), \bF(\lambda)$), where
$$
F(\lambda) = J(\lambda)\Theta \quad \mbox{ and }\quad
\bF(\lambda) = \bTheta J^{-1}(\lambda)
$$
defines a twist for $H=\End_k(A)\otimes_k U$.
Thus, every dynamical twist $J(\lambda)$ for a Hopf algebra $U$
gives rise to a quantum groupoid $H_J = H_{J(\lambda)\Theta}$.
     
\begin{remark}
According to Proposition~\ref{twisting}, the antipode $S_J$ of $H_J$
is given by $S_J(h) = v^{-1}S(h)v$ for all $h\in H_J$, where $S$
is the antipode of $H$ and
$$
v = \sum_{\lambda \mu} E_{\lambda+\mu \lambda}
    (S(J^{(1)}) J^{(2)})(\lambda) P_\mu.
$$
\end{remark}

\begin{remark}
There is a procedure dual to the one described above. It was shown in
\cite{EV} that given a dynamical twist  $J(\lambda)$ it is possible 
to deform the multiplication on the vector space 
$D = \Map(\mathbb{T}\times\mathbb{T},\, k)\otimes_k U^*$ and 
obtain a {\em dynamical quantum group} $D_{J}$. The relation between
$H_J$ and $D_{J}$ in the case when $\dim{U}<\infty$
was established in \cite{EN}, where it was proven that
$D_J$ is isomorphic to $H_J^{*op}$.
\end{remark}

\subsection{Dynamical twists for $U_q(\g)$ at roots of $1$}

Suppose that $\g$ is a simple Lie algebra of type
$A$, $D$ or $E$ and $q$ is a primitive $\ell$th root of unity in $k$,
where $\ell\geq 3$ is odd and coprime with the determinant
of the Cartan matrix $(a_{ij})_{ij=1,\dots,m}$ of $\g$.

Let $U=U_q(\g)$ be the corresponding quantum group
which is a finite dimensional Hopf algebra with generators
$E_i,F_i, K_i$, where $i=1,\dots,m$ and relations as in \cite{L}.

Let $\mathbb{T} \cong (\mathbb{Z}/\ell\mathbb{Z})^m$ be the abelian
group generated by $K_i,\,i=1,\dots,m$. For any $m$-tuple of integers
$\lambda =(\lambda_1,\dots,\lambda_m)$
we will write $K_\lambda = K_1^{\lambda_1}\dots K_m^{\lambda_m}\in\mathbb{T}$.

Let $\R$ be the universal $R$-matrix of $U_q(\g)$ and $\Omega$
be the ``Cartan part'' of $\R$ (\cite{T}).

For arbitrary non-zero constants $\Lambda_1,\dots, \Lambda_m$
define a Hopf algebra automorphism $\Lambda$ of $U$ by
setting
\begin{equation*}
\Lambda(E_i) = \Lambda_i E_i,\quad
\Lambda(F_i) = \Lambda_i^{-1} F_i,\quad  \mbox{ and } \quad
\Lambda(K_i) = K_i \quad\mbox{ for all } i=1,\dots,m.
\end{equation*}
We will say that $\Lambda=(\Lambda_1,\dots,\Lambda_m)$ is {\em generic}
if the spectrum of $\Lambda$ does not contain $\ell$th roots of unity.

Note that the algebra $U$ is $\mathbb{Z}$-graded with
\begin{equation}
\deg(E_i) =1, \quad \deg(F_i) = -1, \quad \deg(K_i) = 0,
\qquad i=1,\dots,m,
\end{equation}
and $\deg(XY)= \deg(X)+\deg(Y)$ for all $X$ and $Y$. Of course,
there are only finitely many non-zero components of $U$   
since it is finite dimensional.   

Let $U_+$ be the subalgebra of $U$ generated by
the elements $E_i,\,K_i,\,i=1,\dots,m$, $U_-$ be
the subalgebra  generated by $F_i,\,K_i,\,i=1,\dots,m$, and
$I_\pm$ be the kernels of the projections from $U_\pm$
to the elements of zero degree.

The next Proposition was proven in \cite{ABRR}, \cite{ES}, \cite{ESS}
for generic values of $q$ and in \cite{EN} for $q$ a root of unity :

\begin{proposition}
\label{equation for J}
For every generic $\Lambda$ there exists a unique element
$J(\lambda) \in 1 + I_+ \otimes I_-$ that satisfies the
following ABRR relation \cite{ABRR}, \cite{ES}, \cite{ESS} :
\begin{equation}
\label{eqn: ABRR}
(\Ad K_\lambda \circ \Lambda \otimes \id) (\R J(\lambda) \Omega^{-1})
= J(\lambda),\quad  \lambda\in \mathbb{T}.
\end{equation}
The element
\begin{equation}
\tJ(\lambda)= J(2\lambda+h^{(1)}+h^{(2)})
\end{equation}
is a dynamical twist for $U_q(\g)$
in the sense of Definition~\ref{dynamical twist}.
\end{proposition}
\begin{proof}
The existence and uniqueness of the solution of (\ref{eqn: ABRR})
and the fact that $\tJ(\lambda)$ satisfies  (\ref{dynamical equation})
are established by induction on the degree of the first component
of $J(\lambda)$, see (\cite{EN}, 5.1).
\end{proof}

\begin{remark}
(a) The reason for introducing the ``shift'' automorphism $\Lambda$
is to avoid singularities in equation \cite{ABRR}.  Thus, we have
a family of dynamical deformations of $U_q(\g)$ (and, therefore,
a family of quantum groupoids depending on $m$ parameters). 
\newline (b) One can generalize the above construction and associate
a  dynamical twist $\tJ_T(\lambda)$ with any 
{\em generalized Belavin-Drinfeld triple} which consists of subsets
$\Gamma_1,\Gamma_2$ of the set $\Gamma=(\alpha_1,\dots,\alpha_m)$ 
of simple roots of $\g$ together with an inner product preserving
bijection $T: \Gamma_1 \to \Gamma_2$, see \cite{EN}, \cite{ESS}.
\end{remark}

\begin{example}
\label{sl 2}
Let us give an explicit expression for the twists
$J(\lambda)$ and $\tJ(\lambda)$ in the case
$\g=sl(2)$. $U_q(\g)$ is then generated by $E,F,K$ with the standard
relations. The element analogous to
$J(\lambda)$ for generic $q$ was computed
already in \cite{Ba} (see also \cite{BBB}).
If we switch to our conventions, this element will take the form
$$
J(\lambda)=\sum_{n=0}^\infty q^{-n(n+1)/2}\frac{(1-q^2)^n}{[n]_q!}
(E^n\otimes F^n)
\prod_{\nu=1}^n\frac{ \Lambda q^{2\lambda}}{1-\Lambda q^{2\lambda+2\nu}
(K\otimes K^{-1})}.
$$
It is obvious that the formula for
$q$ being a primitive $\ell$-th root of unity is simply obtained
by truncating this formula:
$$
J(\lambda)=\sum_{n=0}^{\ell-1} q^{-n(n+1)/2}\frac{(1-q^2)^n}{[n]_q!}
(E^n\otimes F^n)
\prod_{\nu=1}^n\frac{ \Lambda q^{2\lambda}}{1-\Lambda q^{2\lambda+2\nu}
(K\otimes K^{-1})}.
$$
Therefore,
$$
\tJ(\lambda)=\sum_{n=0}^{\ell-1} q^{-n(n+1)/2}\frac{(1-q^2)^n}{[n]_q!}
(E^n\otimes F^n)
\prod_{\nu=1}^n\frac{\Lambda q^{4\lambda}K\otimes K}{1-\Lambda q^{4\lambda+2\nu}
(K^2\otimes 1)}.   
$$
\end{example}

An important new property of the resulting quantum groupoids
$H_\tJ = U_q(\g)_\tJ$, compared with quantum groups $U_q(\g)$ is their
selfduality.   
A twisting of a quasitriangular quantum groupoid 
is again quasitriangular (Remark~\ref{rem: twist}(c)), so the 
twisted $R$-matrix
\begin{equation*}
\R(\lambda)
= \bTheta_{21}\R \Theta
= \sum_{\lambda\mu\nu}\, E_{\lambda \lambda+\nu}P_\mu
\R^{\tJ(1)}(\lambda) \otimes E_{\lambda+\mu \lambda} 
\R^{\tJ(2)}(\lambda)P_\nu,
\end{equation*}
where $\R^{J}(\lambda) = \tJ^{-1}_{21}(\lambda)R \tJ(\lambda)$,
establishes a homomorphism between $D_\tJ = H_\tJ^{*op}$ and $H_\tJ$.
One can show that the image of $\R(\lambda)$ contains all the generators
of $U_q(\g)$,  i.e., that the above homomorphism is in fact an isomorphism, 
see (\cite{EN}, 5.3).

\begin{remark}
The same selfduality result holds for any generalized Belavin-Drinfeld triple
for which $T$ is an automorphism of the Dynkin diagram of $\g$ 
(\cite{EN}, 5.4).
\end{remark}

\end{section}


\begin{section}
{Semisimple and $C^*$-quantum groupoids}

\subsection{Definitions}
\begin{definition}
\label{semisimple and C*}
A quantum groupoid is said to be \textit{semisimple} (resp., 
{\em $*$- or $C^*$-quantum groupoid}) if its algebra $H$ is semisimple 
(resp., a $*$-algebra over a field $k$ with involution, or finite-dimensional
$C^*$-algebra over the field $\mathbb{C}$ and $\Delta$ is a $*$-homomorphism). 
\end{definition}
 
Groupoid algebras and their duals give examples of commutative
and cocommutative semisimple quantum groupoids 
(Corollary~\ref{kG is semisimple}),  which are 
$C^*$-quantum groupoids if the ground  field is $\mathbb{C}$ 
(in which case $g^*=g^{-1}$ for all $g\in G$).

We will describe the class of quantum groupoids possessing \textit{Haar
integrals}, i.e., normalized two-sided integrals (note that if such an integral 
exists then it is unique and is an $S$-invariant idempotent). 

\begin{definition} (\cite{W}, 1.2.1).
\label{W-index} 
Given an inclusion of unital $k$-algebras $N\subset M$, a \textit{conditional 
expectation} $E:M\to N$ is an $N-N$ bimodule map $E$ such that $E(1)=1$. 
A \textit{quasi-basis} for $E$ is an element 
$\sum_i x_i\otimes y_i\in M\otimes_k M$ such that  
$$
\sum_i E(mx_i)y_i=m=\sum_i x_iE(y_i m), \qquad \mbox{ for all } m\in M.
$$

One can check that $\Index E=\sum_i x_i y_i\in Z(M)$ does not depend on a 
choice of a quasi-basis; this element is called the \textit{index} of $E$. 

One can apply this definition to a non-degenerate 
functional $f$ on $H$. If $(l,\lambda)$ is a dual pair of left integrals, 
then $l\2\otimes S^{-1}(l\1)$ is a quasi-basis for $\lambda$ and 
$$
\Index\lambda=S^{-1}\cdot\eps_t(l)\in H_s\cap Z(H).
$$
So a non-degenerate left integral is normalized iff its dual has index $1$.
\end{definition}

\subsection{Existence of the Haar integral}

Given a dual bases $\{f_i\}$ and $\{\xi^i\}$ of $H$ and $\widehat H$, 
respectively,  let us consider the following canonical element
$$
\chi =\sum_i(\xi^i\actl S^{-2}(f_i)).
$$
One can show that $\chi$ is a left integral in $\widehat H$ such that
$\chi(xy)=\chi(yS^2(x))$ for all $x,y\in H$ and $l\act\chi={\widehat S}^2(\widehat 1\actl l)$ 
for any left integral $l$ in $H$ (see \cite{BNSz}, 3.25).

\begin{lemma}[\cite{BNSz}, 3.26]
\label{Haar}
\begin{enumerate}
\item[(i)] The Haar integral $h\in H$ exists iff the above mentioned $\chi$ is
non-degenerate, in which case $(h,\chi)$ is a dual pair of left integrals.
\item[(ii)] A left integral $l$ is a Haar integral iff $\eps_s(l)=1$.
\end{enumerate}
\end{lemma}
\begin{proof} (ii) If $\eps_s(l)=1$, then by 
$l\act\chi={\widehat S}^2(\widehat 1\actl l)$ one has $l\act\chi=\widehat 1$. 
Therefore, $(h,\chi)$ is a dual pair of non-degenerate left integrals. 
The property $\chi(xy)=\chi(yS^2(x))$ is equivalent to that the 
quasi-basis of $\chi$ satisfies
$$
l\2\otimes S^{-1}(l\1)=S(l\1)\otimes l\2,
$$
from where $\Delta(l)=\Delta(S(l))$ and $l=S(l)$. Furthermore, $\eps_t(l)=
\eps_t(S(l))=S\cdot \eps_s(l)=1$. Thus, $l$ is a Haar integral. The inverse 
statement is obvious.
\newline
(i) The "only if" part follows from the proof of (ii). 
If $\chi$ is non-degenerate 
and $h$ is its dual left integral, then, as above, $S(h)=h$; so $h$ is a 
two-sided integral and since $l\act\chi={\widehat S}^2(\widehat 1\actl l)$, 
it is normalized.
\end{proof}
 
\begin{theorem}[\cite{BNSz}, 3.27]
\label{Haar theorem}
Let $H$ be a finite quantum groupoid over an algebraically closed field $k$. 
Then the following conditions are equivalent:
\begin{enumerate}
\item[(i)] There exists a Haar integral, 
\item[(ii)] $H$ is semisimple and there exists an invertible element $g\in H$ 
such that $gxg^{-1}=S^2(x)$ for all $x\in H$ and
$Tr(\pi_\alpha(g^{-1}))\neq 0$ for all irreducible representations $\pi_\alpha$
of $H$ (here $Tr$ is the usual trace on a matrix algebra). 
\end{enumerate}
\end{theorem}
\begin{proof} 
The assumption on $k$ is used only to ensure that
$H=\oplus_\alpha M_{n_\alpha}(k)$, once knowing that it is semisimple, 
so that there is a $k$-basis of matrix units $\{e^{\alpha}_{ij}\}$ in $H$.

(ii)$\implies$(i):  Recall that any trace $\tau$ on $H$ is 
completely determined by its \textit{trace vector}  
$\tau_\alpha =\tau(1_\alpha)$. 
Now let $\tau:H\to k$ be the trace with trace vector $\tau_\alpha=
\tr \pi_\alpha(g^{-1})$. Then $\tau$ is non-degenerate and has a quasi-basis
$$
\sum_i x_i\otimes y_i=\sum_\alpha{\frac{1}{\tau_\alpha}}\sum^{n_\alpha}_{i,j=1}
e^{\alpha}_{ij}\otimes e^{\alpha}_{ji}.
$$
Notice that $\sum_i x_ig^{-1}y_i=1$. It is straightforward to verify that
$\chi'=g\act\tau$ coincides with $\chi$, so $\chi$ is  
non-degenerate and therefore its dual left integral $l$ has $\eps_s(l)=1$ by 
$l\act\chi={\widehat S}^2(\widehat 1\actl l)$. Thus, $l$ is a Haar integral.

(i)$\implies$(ii): If $h$ is a Haar integral, then $H$ is semisimple by
Theorem~\ref{maschke}. Let $\tau$ be a non-degenerate trace on $H$, 
then there exists a  unique $i\in H$ such that $i\act\tau=\widehat 1=\tau\actl i$. 
One can verify that $i$  is a two-sided non-degenerate integral in $H$ 
and that $S(i)=i$ (see \cite{BNSz}, I.3.21). 
To prove that $S^2$ is inner, it is enough to construct a non-degenerate 
functional $\chi$ on $H$ such that $\chi(xy)=\chi(yS^2(x))$ for all $x,y\in H$. 
But the proof of Lemma \ref{Haar} shows that this is the case for the dual 
left  integral to $i$.
\end{proof}

\begin{remark} 
\label {C*-semisimple}
\begin{enumerate}
\item[(a)]
There is a unique normalization of $g$
from Theorem~\ref{Haar theorem}(ii) 
such that the following conditions hold (\cite{BNSz}, 4.4) :
\begin{enumerate}
\item[(i)] $\tr(\pi_\alpha(g^{-1}))=\tr (\pi_\alpha(g))$ for all irreducible
representations $\pi_\alpha$ of $H$ (here $\tr$ is a usual trace on a
matrix algebra);
\item[(ii)] $S(g)=g^{-1}$.
\end{enumerate}
One can show that such a $g$ is {\em group-like}, i.e.,
$$
\Delta(g) =(g\otimes g)\Delta(1) = \Delta(1) (g\otimes g).
$$
The element $g$ implementing $S^2$ and satisfying
normalization conditions (i) and (ii) is called
the {\em canonical group-like element} of $H$.
\item[(b)]
Since $\chi = g\act\tau$ is the dual left integral to $h$, its quasi-basis
$\sum_i x_ig^{-1}\otimes y_i$ equals to $h\2\otimes S^{-1}(h\1)$, which
implies the formula
$$ 
(S\otimes id)\Delta(h)=\sum_i x_i\otimes g^{-1}y_i=\sum_\alpha 1/\tau_\alpha
\sum_{ij}e^\alpha_{ij}g^{-1/2}\otimes g^{-1/2}e^\alpha_{ji}.
$$
\end{enumerate}
\end{remark}

A dual notion to that of the Haar integral is the Haar measure. 
Namely, a functional $\phi \in \widehat H$ is said to be a {\em Haar measure} 
on $H$  if it is a normalized left and right invariant measure and 
$\phi\circ S=\phi$. Any of the  
equivalent conditions of Theorem \ref{Haar} is also equivalent to existence 
(and uniqueness) of the Haar measure.

\subsection{$C^*$-quantum groupoids}

Definition \ref{semisimple and C*} and the uniqueness of the unit, counit and 
the antipode (see Proposition \ref{properties of S}) imply that 
$$
1^*=1,\quad \eps(x^*)= \overline{\eps(x)},\quad (S\circ *)^2 =\id
$$ 
for all $x$ in any $*$-quantum groupoid. It is also easy to check the relations 
$$
\eps_t(x)^*=\eps_t(S(x)^*),\qquad \eps_t(x)^*=\eps_t(S(x)^*),
$$
therefore, $H_t$ and $H_s$ are $*$-subalgebras, and to show that the dual, 
$\widehat H$, is also a $*$-quantum groupoid with respect to the $*$-operation
\begin{equation}
\la\phi^*,x\ra=\overline{\la\phi,S(x)^*\ra}\qquad \mbox{ for all } 
\phi\in \widehat H,\ x\in H.
\end{equation}

The $*$-operation allows to simplify the axioms of a quantum groupoid 
(cf. the  axioms used in \cite{NV1}, \cite{N1}). The second parts of
equalities~(\ref{eps m}) and (\ref{Delta 1})  
of Definition~\ref{finite quantum groupoid} 
follow from the rest of the axioms, also $S*\id = \eps_s$ is 
equivalent to $\id* S=\eps_t$. Alternatively, 
under the condition that the antipode is both algebra and coalgebra 
anti-homomorphism, the axioms (\ref{eps m}) and (\ref{Delta 1}) can be
replaced  by the identities of Proposition~\ref{counital properties} (ii) 
and (v) involving the target counital map.

\begin{theorem} (\cite{BNSz}, 4.5)
\label{C*-Haar}
In a $C^*$-quantum groupoid the Haar integral $h$ exists, $h=h^*$ and
$$
(\phi,\psi)=\la\phi^*\psi,h\ra,\qquad \phi,\psi\in \widehat H
$$
is a scalar product making $\widehat H$ a Hilbert space where the left regular
representation of $\widehat H$ is faithful. 
Thus, $\widehat H$ is a $C^*$-quantum groupoid, too.
\end{theorem}
\begin{proof}
Clearly, $H$ and $g$ verify all 
the conditions of Theorem \ref{Haar theorem}, from where the existence of Haar 
integral follows. Since $h$ is non-degenerate, the scalar product 
$(\cdot,\cdot)$ is  also non-degenerate. By the equality 
$$
(\phi,\phi)=\la\phi^*\phi,h\ra=\overline{\la\phi,S(h\1)^*\ra}\la\phi,h\2\ra,
$$
positivity of $(\cdot,\cdot)$ follows from Remark \ref{C*-semisimple}(b).
\end{proof}

We will denote by $\widehat h$ the Haar measure of $\widehat H$.

\begin{remark} 
\label{eps is positive}
$\eps$ is a positive functional, i.e., $\eps(x^*x) \geq 0$ for all $x\in H$.
Indeed, for all $x\in H$ we have
$$
\eps(x^*x)
= \eps(x^*1\1)\eps(1\2 1'\2)\eps(1'\1 x)=\eps(\eps_t(x)^* \eps_t(x))
= \la\widehat h,\eps_t(x)^*\eps_t(x)\ra\geq 0,
$$
where $\widehat h\vert_{H_t}=\eps\vert_{H_t}$, since
$\la\widehat h,z\ra=\la\widehat\eps_t(\widehat h),z\ra=\la\widehat 1,z\ra$ 
for all $z\in H_t$.
\end{remark}

The Haar measure provides target and source {\em Haar conditional 
expectations} (all properties are easy to verify):
\begin{eqnarray*}
E_t &:& H\to H_t :  E_t(x)=(id\otimes\widehat h)\Delta(x), \\
E_s &:& H\to H_s :  E_s(x)=(\widehat h\otimes id)\Delta(x).
\end{eqnarray*}
Let us introduce an element $g_t=E_t(h)^{1/2}\in H_t$.

\begin{remark}
\label{group-like}
(\cite{BNSz}, 4.12).
$g_t$ is positive and invertible, and the canonical group-like element 
of $H$ can be written as $g=g_tS(g_t^{-1})$.
\end{remark}

One can check that for a quasitriangular $*$-quantum groupoid $\bR=\R^*$.
\begin{proposition}[\cite{NTV}, 9.3]
\label{*-double} 
If $H$ is a $C^*$-quantum groupoid, then $D(H)$ is a quasitriangular
$C^*$-quantum groupoid.
\end{proposition}
\begin{proof}
The result follows from the fact that $\widehat{D(H)}$ is a $C^*$-quantum
groupoid, which is easy to verify using the explicit formulas from
Remark~\ref{dual double}.
\end{proof}

\begin{proposition}
\label{*-Ribbon WHA}
A quasitriangular $C^*$-quantum groupoid $H$ is automatically ribbon 
with the ribbon element $\nu=u^{-1}g=gu^{-1}$, where $u$ is the Drinfeld
element from Definition~\ref{Drinfeld element} and $g$ is the canonical
group-like element.
\end{proposition}
\begin{proof}
Since $u$ also implements $S^2$ (Proposition \ref{elements u and v}), 
$\nu=u^{-1}g$ is central, therefore $S(\nu)$ is also central.
Clearly, $u$ must commute with $g$. The same Proposition gives
$\Delta(u^{-1}) = \R_{21}\R(u^{-1}\otimes u^{-1})$,
which allows us to compute
$$
\Delta(\nu) = \Delta(u^{-1})(g\otimes g) = 
\R_{21}\R(u^{-1}g\otimes u^{-1}g) = \R_{21}\R(\nu \otimes \nu).
$$
Propositions~\ref{properties of R} and \ref{elements u and v}
and the trace property imply that
\begin{eqnarray*}
\tr(\pi_\alpha(u^{-1})) 
&=& \tr(\pi_\alpha( \R^{(2)}S^2(\R^{(1)}) )) \\ 
&=& \tr (\pi_\alpha(S^3(\R^{(1)}) S(\R^{(2)}) ) = \tr(\pi_\alpha(S(u^{-1}))).
\end{eqnarray*}
Since $u^{-1} = \nu g^{-1}$ and $\nu$ is central, the above relation 
means that 
$$
\tr(\pi_\alpha(\nu))\tr(\pi_\alpha(g^{-1})) =
\tr(\pi_\alpha(S(\nu)))\tr(\pi_\alpha(g)),
$$
and, therefore, $\tr(\pi_\alpha(\nu)) = \tr(\pi_\alpha(S(\nu)))$ for
any irreducible representation $\pi_\alpha$, which shows that
that $\nu = S(\nu)$ (cf. \cite{EG}).
\end{proof}

\begin{remark}
\label{C*-q-trace}
For a connected ribbon $C^*$-quantum groupoid we have:
$$
\tr_q(f) = (\dim H_t)^{-1} \Tr_V(g\circ f),\ \dim_q(V) = 
(\dim H_t)^{-1} \Tr_V(g).
$$
for any $f\in \End(V)$, $V$ is an $H$-module. 
\end{remark}

\subsection{$C^*$-quantum groupoids and unitary modular categories}

To define the (unitary) representation category $\URep(H)$ of a $C^*$-quantum groupoid
$H$ we consider {\it unitary} $H$-modules, i.e., $H$-modules $V$ 
equipped with a scalar product 
$$
(\cdot,\cdot): V\times V\to \mathbb{C} \qquad
\mbox{ such that } \qquad 
(h\cdot v,w)=(v,h^*\cdot w)\ \forall h\in H, v,w\in V.
$$
The notion of a morphism in this category is the same as in $\Rep(H)$.

The monoidal product of $V,W\in \URep(H)$ is defined as follows. We
construct a tensor product $V\otimes_\mathbb{C} W$ of Hilbert spaces and note
that the action of $\Delta(1)$ on this left $H$-module is an orthogonal
projection. The image of this projection is, by definition, the
monoidal product 
of $V,W$ in $\URep(H)$. Clearly, this definition is compatible with the monoidal
product of morphisms in $\Rep(H)$.
 
For any $V\in \URep(H)$, the dual space $V^*$ is naturally identified 
($v\to \overline v$) with the conjugate Hilbert space, and under this 
identification we have 
$h\cdot\overline v= \overline{S(h)^*\cdot v}\ (v\in V, \overline v\in V^*$). 
In this way $V^*$ becomes a unitary $H$-module with scalar product 
$(\overline v, \overline w)=(w,gv)$, where $g$ is the canonical group-like 
element of $H$. 

The unit object in $\URep(H)$ is $H_t$ equipped with scalar product 
$(z,t)_{H_t}= \eps(zt^*)$ (it is known \cite{BNSz}, \cite{NV1} that
the restriction of $\eps$ to $H_t$ is a non-degenerate positive form). 
One can verify that the maps  $l_V,r_V$ and their inverses are isometries
and that $\URep(H)$ is a monoidal category with duality
(see also \cite{BNSz}, Section 3). In a natural way we have a notion of a 
{\em conjugation} of morphisms (\cite{T}, II.5.1).

\begin{remark}
\label{*-R-matrix}
a) One can check that for a quasitriangular $*$-quantum groupoid $H$
the braiding is an isometry in $\URep(H) : c^{-1}_{V,W}=c^{*}_{V,W}$.

b) For a ribbon $C^*$-quantum groupoid $H$,
the twist is an isometry in $\URep(H)$. 
Indeed, the relation $\theta_{V}^*=\theta_{V}^{-1}$ 
is equivalent to the identity $S(u^{-1}) = u^*$, 
which follows from Proposition~\ref{properties of R} and 
Remark~\ref{*-R-matrix}a).
\end{remark} 

A {\em Hermitian ribbon category} is a ribbon category endowed with a 
conjugation of morphisms $f\mapsto\overline{f}$ satisfying natural conditions
of (\cite{T}, II.5.2).
A {\em unitary ribbon category} is a Hermitian ribbon category over the
field $\mathbb{C}$ such that for any morphism $f$ we have $\tr_q(f\overline{f})
\geq 0$.
In a natural way we have a {\em conjugation}
of morphisms in $\URep(H)$. Namely, for any morphism $f:V\to W$ we define 
$\overline f:W\to V$ as $\overline f (w)=\overline{f^*(\overline w)}$ for
any $w\in W$. Here $\overline w\in W^*, f^*:W^*\to V^*$ is the standard dual
of $f$ (see \cite{T}, I.1.8) and $\overline{f^*(\overline w)}\in V$.
For the proof of the following lemmma see (\cite{NTV}, 9.7).

\begin{lemma}
\label{hermit-ribbon}
Given a quasitriangular $C^*$-quantum groupoid $H,\ \URep(H)$ is a unitary 
ribbon $Ab$-category with respect to the above conjugation of morphisms.
\end{lemma}

The next proposition extends (\cite{EG}, 1.2).

\begin{proposition}
\label{rep of C* is modular}
If $H$ is a connected $C^*$-quantum groupoid, then $\Rep(D(H))$ is a unitary 
modular category.
\end{proposition}
\begin{proof}
The proof follows from Lemmas~\ref{hermit-ribbon}, 
\ref{factorizable implies modular} and
Propositions~\ref{D(H) is factorizable}, \ref{*-double}.
\end{proof}
\end{section}


\begin{section} 
{$C^*$-quantum groupoids and subfactors : the depth 2 case}
 
In this section we characterize finite index  type II$_1$ subfactors 
of depth $2$ in terms of $C^*$-quantum groupoids.
Recall that a II$_1$ factor is a $*$-algebra of operators on a Hilbert
space that coincides with its second centralizer, or, equivalently,
is weakly closed (i.e., a von Neumann algebra)
with the trivial center that admits a finite trace.
An  example of such a factor is the group von Neumann algebra of a discrete
group whose every non-trivial conjugacy class is infinite, the corresponding
trace given by its Haar measure. It is also known, that there exists
a unique, up to an isomorphism, hyperfinite (i.e., generated by an 
increasing sequence of finite-dimensional $C^*$-algebras) II$_1$ factor.

There is a notion of index for subfactors, extending the notion of index of
a subgroup in a group, though it can be non-integer. 
For the foundations of  the subfactor theory see \cite{GHJ}, \cite{JS}.

\subsection{Actions of $C^*$-quantum groupoids on von Neumann algebras}
Let a von Neumann algebra $M$ be a left 
$H$-module algebra in the sense of Definition \ref{module algebra} via a
weakly continuous action of a $C^*$-quantum groupoid 
\begin{equation*}
H\otimes_\mathbb{C} M \ni  x\otimes m \mapsto (x\lact m)\in M
\end{equation*}
such that $(x\lact m)^* = S(x)^* \lact m^*$ 
and $x\lact 1 =0$ iff $\eps_t(t)=0.$

Then one can show (\cite{NSzW}, 3.4.2) that the smash product algebra 
(now we call it {\em crossed product algebra} and  denote by $M\rtimes H$),  
equipped with an involution  
$[m\otimes x]^* = [(x\1^*\lact m^*) \otimes  x\2^*]$, is a von Neumann algebra.

The collection
$ M^H =\{ m\in M \mid x\lact m = \eps_t(x)\lact m,\ \forall x\in H\}$
is a von Neumann subalgebra of $M$, called a {\em fixed point subalgebra}.
The centralizer $M^\prime \cap M \rtimes H$ always contains
a the source counital subalgebra $H_s$. 
Indeed, if $y\in H_s$, then $\Delta(y) =1\1\otimes 1\2 y$, therefore
\begin{eqnarray*}
[1\otimes y][m\otimes 1]
&=& [(y\1\lact m)\otimes y\2]=[(1\1\lact m)\otimes  1\2 y] \\
&=& [m\otimes y]=[m\otimes 1][1\otimes y],
\end{eqnarray*}
for any $m\in M$, and $H_s\subset M^\prime \cap M \rtimes H$.
An action of $H$ is
called {\em minimal} if $H_s = M^\prime \cap M\rtimes H$.

One can define the dual action of $\widehat H$  on the von
Neumann algebra $M\rtimes H$ as in (\ref{dual action}) and 
construct the von Neumann algebra $(M\rtimes H) \rtimes \widehat H$.

\begin{remark}
The following results hold true if $S^2=\id$, i.e., if $H$ is involutive
(\cite{N1}, 4.6, 4.7).
\begin{enumerate}
\item[(a)] $M\subset M\rtimes H\subset (M\rtimes H)
\rtimes \widehat H$ is a basic construction;
\item[(b)] $M\rtimes H$ is free as a left $M$-module;
\item[(c)]  $(M\rtimes H)\rtimes \widehat H \cong M\otimes_\mathbb{C}
M_n(\mathbb{C})$ for some integer $n$.
\end{enumerate}
\end{remark}

\begin{example} 
\label{Weyl algebra} 
For a trivial action of $H$ on $H_t$  and
the corresponding dual action of $\widehat H$ on $H\cong H_t\rtimes H$ (see
Example \ref{examples of actions} (i) and (ii)), 
$H_t\subset H\subset H\rtimes \widehat H$ is the basic
construction of finite-dimensional $C^*$-algebras with respect to the
Haar conditional expectation $E_t$ (\cite{BSz2}, 4.2). 
\end{example}

If $H$ is a connected
$C^*$-quantum groupoid having a minimal action on a II${}_1$ factor $M$
then one can show  (\cite{NSzW}, 4.2.5, 4.3.5) that 
$\widehat H$ is also connected and that 
$$ 
N=M^H\subset M\subset M_1=M \rtimes H\subset M_2=(M\rtimes H)\rtimes \widehat H
\subset \cdots
$$
is the Jones tower of factors of finite index with {\em the derived tower} 
$$
N'\cap N=\mathbb{C}\subset N'\cap M=H_t\subset N'\cap M_1=H\subset N'\cap M_2=
H\rtimes \widehat H\subset  \cdots
$$ 
The fact that the last triple of finite-dimensional $C^*$-algebras is a
basic construction, means exactly that the subfactor $N\subset M$ 
has depth 2. Moreover, the finite-dimensional $C^*$-algebras 
\begin{equation}
\label{square}
\begin{array}{ccc}
H^* & \subset & H^* \rtimes H  \\
\cup&         & \cup        \\   
H_t= H_s^* & \subset & H,
\end{array}
\end{equation}
form {\em a canonical commuting square}, 
which completely  determines the equivalence class of the
initial subfactor. This implies that any biconnected $C^*$-quantum groupoid
has at most one minimal action on a given II$_1$ factor and thus corresponds 
to no more than one (up to equivalence) finite index depth 2 subfactor.

\begin{remark}
It was shown in (\cite{N1}, 5) that any biconnected involutive $C^*$-quantum 
groupoid has a minimal action on the  hyperfinite II$_1$ factor. This
action is constructed by iterating the basic construction for the
square (\ref{square}) in the horizontal direction, see \cite{N1} for
details.
\end{remark}

\subsection{Construction of a $C^*$-quantum groupoid from a depth 2 subfactor}

Let $N\subset M$ be a finite index ($[M:N]= \lambda^{-1}$) depth 2
II${}_1$ subfactor and
\begin{equation}
N \subset M \subset M_1 \subset M_2 \subset \cdots
\end{equation}
the corresponding Jones tower, $M_1 = \langle M,\,e_1 \rangle,\,
M_2 = \langle M_1,\,e_2 \rangle,\dots,$ where $e_1 \in N^\prime
\cap M_1,\,e_2 \in M^\prime \cap M_2,\cdots$
are the Jones  projections. The depth 2 condition means that
$N^\prime \cap M_{2}$ is the basic construction of the inclusion
$N^\prime \cap M  \subset  N^\prime \cap M_{1}$. Let $\tau$ be
the trace on $M_2$ normalized by $\tau(1)=1$.

Let us denote
$$
A=N'\cap M_1, \quad B=M'\cap M_2, \quad B_t=M'\cap M_1, \quad B_s=M_1'\cap M_2
$$
and let $\Tr$ be the trace of the regular representation of $B_t$ on
itself. Since both $\tau$ and $\Tr$ are non-degenerate, there exists a positive
invertible element $w\in Z(B_t)$  such that $\tau(wz)=\Tr(z)$ for all
$z\in B_t$ (the index of $\tau\vert_{M'\cap M_1}$
in the sense of \cite{W}).

\begin{proposition} 
\label{duality}
There is a canonical non-degenerate duality form between 
$A$ and $B$ defined by 
\begin{equation}
\la a,\, b \ra = \lambda^{-2}\tau(ae_2e_1wb)= \lambda^{-2}Tr E_{B_t}(bae_2e_1),
\end{equation}
for all $a\in A$ and $b\in B$.
\end{proposition}
\begin{proof}
If $a\in A$ is such that $\la a,\, B\ra = 0$ , then, using 
\cite{PP1}, 1.2, one has
$$
\tau(a e_2 e_1 B) = \tau(a e_2 e_1 (N^\prime \cap M_2)) = 0,
$$
therefore, using the properties of $\tau$ and Jones projections, we get
\begin{eqnarray*}
\tau(a a^*) &=& \lambda^{-1} \tau(a e_2 a^*) =
\lambda^{-2}\tau(a e_2 e_1 (e_2 a^* ) ) = 0,
\end{eqnarray*}
so $a=0$. Similarly for $b\in B$.
\end{proof}

Using the duality $\la a,\, b \ra$, one defines a coalgebra structure on $B$ :
\begin{eqnarray}
\la a_1\otimes a_2,\, \Delta(b) \ra &=&   \la a_1 a_2,\, b \ra, \\
\eps(b) &=& \la 1,\, b\ra, 
\end{eqnarray}
for all $a,a_1,a_2\in A$ and $b\in B$. Similarly, one defines a
coalgebra structure on $A$.

We define a linear endomorphism $S_B:B\to B$ by 
\begin{equation}
E_{M_1}( b e_1 e_2) = E_{M_1}(e_2 e_1 S_B(b)),
\end{equation}
Note that $\tau\circ S_B =\tau$.
Similarly one can define a linear endomorphism $S_A:A\to A$ such that
$E_{M^\prime}(S_A(a)e_2 e_1) = E_{M^\prime}(e_1 e_2 a)$
and $\tau\circ S_A =\tau$.

\begin{proposition}
\label{properties of S (1)} (\cite{NV2}, 4.5, 5.2):
The following identities hold :
\begin{enumerate}
\item[(i)] $S_B(B_s)=B_t$,
\item[(ii)] $S_B^2(b) =b$ and $S_B(b)^* = S_B(b^*)$,
\item[(iii)] $\Delta(b)(S_B(z)\otimes 1) = \Delta(b)(1\otimes z), \qquad 
\Delta(bz)= \Delta(b)(z\otimes 1),\ (\forall z\in B_t)$,
\item[(iv)] $S_B(bc) = S_B(c)S_B(b)$, \newline
      $\Delta(S_B(b)) = (1\otimes w^{-1})(S_B(b\2)\otimes S_B(b\1))
      (S_B(w)\otimes 1)$.
\end{enumerate}
\end{proposition}

Define an  antipode and new involution of $B$ by
\begin{equation}
S(b)=S_B(wbw^{-1}), \qquad
b^\star = S(w)^{-1}b^*S(w).
\end{equation} 

\begin{theorem}
\label{B is a quantum groupoid}
With the above operations $B$ becomes a biconnected $C^*$-quantum groupoid 
with the counital subalgebras $B_s$ and $B_t$.  
\end{theorem}
\begin{proof}
One can check (\cite{NV2}, 4.6) that $\Delta(1)\in B_s\otimes B_t$, hence, 
since  $\Delta(1)\in B_s\otimes B_t$ and $B_t$ commutes with $B_s$, we have
\begin{eqnarray*}
(\id\otimes\Delta)\Delta(1) 
&=& 1\1 \otimes \Delta(1)(1\2 \otimes 1) \\
&=& (1 \otimes \Delta(1))(\Delta(1)\otimes 1)
    = (\Delta(1)\otimes 1) (1 \otimes \Delta(1))
\end{eqnarray*}
which is axiom (\ref{Delta 1}) of Definition~\ref{finite quantum groupoid}.
The axiom~(\ref{eps m})  dual to it is obtained similarly, considering
the comultiplication in $A$.

As a consequence of the ``symmetric square'' relations $BM_1 =M_1B =M_2$
one obtains (\cite{NV2}, 4.12) the identity
\begin{equation}
\label{crucial}
bx =\lambda^{-1} E_{M_1}(b\1 x e_2) b\2 \quad \mbox{ for all } \quad
b\in B,\,  x\in M_1,
\end{equation}
from where it follows that
\begin{equation}
\label{good property}
E_{M_1}(bxye_2)  =\lambda^{-1} E_{M_1}(b\1 x e_2) E_{M_1}(b\2ye_2) 
\quad \mbox{ for all } \quad b\in B,\,x,y\in M_1.
\end{equation}
We use this identity to prove that $\Delta$ is a homomorphism :
\begin{eqnarray*}
\la a_1a_2,\, bc \ra 
&=& \la \lambda^{-1}E_{M_1}(ca_1a_2e_2),\, b \ra \\
&=& \la \lambda^{-2}E_{M_1}(c\1 a_1e_2)E_{M_1}(c\2 a_2e_2),\, b \ra \\
&=& \la \lambda^{-1}E_{M_1}(c\1 a_1e_2), b\1 \ra
    \la \lambda^{-1}E_{M_1}(c\2 a_2 e_2),\, b\2 \ra    \\
&=& \la a_1,\, b\1c\1 \ra \la a_2,\, b\2 c\2 \ra,
\end{eqnarray*}
for all $a_1,a_2 \in A$, whence $\Delta(bc) =  \Delta(b)\Delta(c)$.
Next, 
\begin{eqnarray*}
\la a,\, \eps(1\1b)1\2 \ra
&=& \la 1,\,1\1b \ra \la a,\, 1\2 \ra \\
&=& \la \lambda^{-1}E_{M_1}(be_2), 1\1 \ra \la  a,\, 1\2 \ra \\
&=& \la \lambda^{-1}E_{M_1}(be_2)a, 1 \ra
 = \la a,\, \lambda^{-1}E_{M_1}(be_2) \ra, 
\end{eqnarray*}
therefore $\eps_t(b) = \lambda^{-1}E_{M_1}(be_2)$ (similarly,
$\eps_s(b) = \lambda^{-1}E_{M_1'}(e_2b)$).
Using Equation (\ref{good property}) we have
\begin{eqnarray*}
\la a,\, b\1 S(b\2) \ra 
&=& \la a,\, b\1 S_B(wb\2 w^{-1}) \ra \\
&=&  \lambda^{-3} \tau (E_{M_1}(S_B(w b\2 w^{-1}) a e_2) e_2 e_1 wb\1) \\
&=&  \lambda^{-3} \tau (E_{M_1}(e_2 a wb\2 w^{-1}) e_2 e_1 wb\1) \\
&=& \lambda^{-3} \tau (E_{M_1}(e_2 a wb\2 w^{-1})E_{M_1}(e_2 e_1 wb\1) ) \\
&=& \lambda^{-2} \tau (E_{M_1}(e_2 a e_1 wb)) = \la a,\, \eps_t(b) \ra. 
\end{eqnarray*}
The antipode $S$ is anti-multiplicative and anti-comultiplicative,
therefore axiom~(\ref{S id S}) of Definition~\ref{finite quantum groupoid}
follows from the other axioms. Clearly,
$\eps_t(B) = M'\cap M_1$ and $\eps_s(B) = M_1'\cap M_2$.
$B$ is biconnected, since the inclusion
$B_t=M^\prime \cap M_1 \subset B =M^\prime \cap M_2$  is connected
(\cite{GHJ}, 4.6.3) and $B_t \cap B_s =(M^\prime \cap M_1)\cap
(M_1^\prime \cap M_2)= \mathbb{C}$.

Finally, from the properties of
$\Delta$ and $S_B$ we have $\Delta(b^\star)=\Delta(b)^{\star\otimes\star}$.
\end{proof}

\begin{remark}
\label{properties of B}
\begin{enumerate}
\item[(i)] $S^2(b) = gbg^{-1}$, where $g=S(w)^{-1}w$. 
\item[(ii)] The Haar projection of $B$ is $e_2w$ 
and the normalized Haar functional is $\phi(b) =$ $= \tau(S(w)wb)$.
\item[(iii)]  The non-degenerate duality $\la\, ,\,\ra$ makes 
$A=N^\prime \cap M_1$ the $C^*$-quantum groupoid dual to $B$.
\item[(iv)] Since the Markov trace (\cite{GHJ}, 3.2.4) of the inclusion
$N\subset M$ is also the Markov trace of the finite-dimensional inclusion
$A_t\subset A$, it is clear that $[M:N]=[A:A_t]=[B:B_t]$. We call last
number {\em the index of the quantum groupoid} $B$.

\item[(v)] Let us mention two classification results obtained in 
(\cite{NV2}, 4.18, 4.19)  by quantum groupoid methods:
\subitem (a) If $N\subset M$ is a depth 2 subfactor such that
$[M:N]$ is a square free integer (i.e., $[M:N]$ is an integer which
has no divisors of the form $n^2,\, n > 1$), then 
$N^\prime\cap M=\mathbb{C}$, and there is a (canonical) minimal action
of a Kac algebra $B$ on $M_1$ such that $M_2 \cong  M_1\rtimes B$ 
and $M=M_1^B$.
\subitem (b) If $N\subset M$ is a depth 2  II${}_1~$subfactor such that
$[M:N]=p$ is prime, then $N^\prime\cap M=\mathbb{C}$, and there is an
outer action of the cyclic group $G=\mathbb{Z}/p\mathbb{Z}$ on $M_1$ such that 
$M_2 \cong  M_1\rtimes G$ and $M=M_1^G$.
\end{enumerate}
\end{remark}

\subsection{Action on a subfactor}
Equation~(\ref{good property}) suggests the following definition 
of the action of $B$ on $M_1$.

\begin{proposition}
\label{the action}
The map $\lact : B \otimes_\mathbb{C} M_1 \to M_1$ :
\begin{equation}
b \lact x = \lambda^{-1} E_{M_1}(b x e_2)
\end{equation}
defines a left action of $B$ on $M_1$ (cf.\ \cite{S}, Proposition 17).
\end{proposition}
\begin{proof}
The above map defines a left $B$-module structure
on $M_1$, since $1\lact x =x$ and
$$
b \lact (c \lact x) 
= \lambda^{-2}  E_{M_1}(b  E_{M_1}(c x e_2)  e_2) 
= \lambda^{-1}  E_{M_1}( bc x e_2 ) = (bc) \lact x.
$$
Next, using equation~(\ref{good property}) we get
\begin{eqnarray*}
b \lact xy 
&=& \lambda^{-1} E_{M_1}( b xy e_2)  
= \lambda^{-2} E_{M_1}( b\1 x e_2) E_{M_1}( b\2 y e_2) \\
&=& ( b\1 \lact x)( b\2 \lact y).
\end{eqnarray*}
By \cite{NV2}, 4.4 and properties of $S_B$ we also get
\begin{eqnarray*}
S(b)^\star \lact x^*
&=& \lambda^{-1} E_{M_1}( S(b)^\star x^* e_2 ) 
 = \lambda^{-1} E_{M_1}( S_B(w)^{-1} S_B(wbw^{-1})^* S_B(w) x^* e_2 ) \\
&=& \lambda^{-1} E_{M_1}( S_B(b^*)  x^* e_2 ) 
 =  \lambda^{-1} E_{M_1}( e_2  x^* b^* )  \\
&=&  \lambda^{-1} E_{M_1}( bx e_2)^* = (b \lact x)^*.
\end{eqnarray*}
Finally,
$$
b \lact 1 = \lambda^{-1} E_{M_1}( b e_2) = 
\lambda^{-1} E_{M_1}( \lambda^{-1} E_{M_1}(b e_2)e_2)
= \eps_t(b) \lact 1,
$$ 
and $b\lact 1 =0$  iff $\eps_t(b) = \lambda^{-1} E_{M_1}(b e_2)= 0$.
\end{proof}

\begin{proposition}
\label{fixed point}
$M_1^B = M$, i.e. $M$ is the fixed point subalgebra of $M_1$.
\end{proposition}
\begin{proof}
If $x\in M_1$ is such that $b\lact x = \eps_t(b)\lact x$ for all
$b \in B$, then $ E_{M_1}(bxe_2) = E_{M_1}(\eps_t(b)xe_2) =
E_{M_1}(be_2) x$. Taking $b =e_2$, we get $E_{M}(x) =x$ which
means that $x \in M$. Thus, $M_1^B \subset M$.

Conversely, if $x\in M$, then $x$ commutes with $e_2$ and 
$$
b\lact x = \lambda^{-1} E_{M_1}(b e_2 x) = 
\lambda^{-1} E_{M_1}( \lambda^{-1} E_{M_1}(be_2)e_2 x)
= \eps_t(b) \lact x,
$$
therefore $M_1^B = M$. 
\end{proof} 

\begin{proposition}
\label{crossproduct}
The map $\theta : [x\otimes b] \mapsto x S(w)^\half b S(w)^{-\half}$
defines a von Neumann algebra isomorphism between $M_1\rtimes B$ and $M_2$.
\end{proposition}
\begin{proof}
By definition of the action $\lact$ we have :
\begin{eqnarray*}
\theta([x(z\lact 1) \otimes b]) 
&=& x S(w)^{\half} \lambda^{-1} E_{M_1}(ze_2)b S(w)^{-\half} \\
&=& x S(w)^{\half} zb S(w)^{-\half} = \theta([x\otimes zb]),
\end{eqnarray*}
for all $x \in M_1,\, b\in B,\, z\in B_t$, so
$\theta$ is a well defined linear map from
$M_1\rtimes B = M_1 \otimes_{B_t} B$  to $M_2$.
 It is surjective since  an orthonormal basis of 
 $B=M^\prime \cap M_2$
over $B_t=M^\prime \cap M_1$ is also a basis of $M_2$ over $M_1$
(\cite{Po}, 2.1.3).
Finally, one can check that 
$\theta$ is a von Neumann algebra isomorphism (\cite{NV2}, 6.3).
\end{proof}
  
\begin{remark}
\label{minimality}
\begin{enumerate}
\item[(i)]
The action of $B$ constructed in Proposition~\ref{the action}
is minimal, since we have $M_1^\prime \cap M_1\rtimes B = M_1^\prime \cap M_2 
=B_s$ by Proposition~\ref{crossproduct}. 
\item[(ii)]
If $N'\cap M=\mathbb{C}$, then $B$ is a usual
Kac algebra (i.e., a Hopf $C^*$-algebra) and we recover the well-known
result proved in \cite{S}, \cite{L}, and \cite{D}.
\item [(iii)]
It is known (\cite{Po}, Section 5) that 
there are exactly 3 non-isomorphic subfactors of depth 2 
and index 4 of the hyperfinite II${}_1$ factor $R$ such that 
$N'\cap M \neq \mathbb{C}$. 
One of them is $R\subset R\otimes_\mathbb{C} M_2(\mathbb{C})$, 
two others can be viewed as diagonal subfactors of the form 
$$
\{  \begin{pmatrix} x & 0 \cr 0 & \alpha(x) \end{pmatrix}
\vert x\in R \} \subset R\otimes M_2(\mathbb{C}),
$$
where $\alpha$ is an outer automorphism of $R$  such that  
$\alpha^2 =\id$, resp.\ $\alpha^2$ is inner and $\alpha^2 \neq\id$. 

An explicit description of the corresponding quantum groupoids
can be found in (\cite{NV1}, 3.1, 3.2), (\cite {NV2}, 7). 
Also, in (\cite {NV2}, 7)  we describe the $C^*$-quantum groupoid 
$B=M_2(\mathbb{C}) \oplus M_3(\mathbb{C})$
corresponding to the subfactor with the principal graph $A_3$ (\cite{GHJ}, 
4.6.5) having  the index $\lambda^{-1} = 16\cos^4\frac{\pi}{5}$.
\end{enumerate}
\end{remark}

\end{section}


\begin{section} 
{$C^*$-quantum groupoids and subfactors : the finite depth case}

Here we show that
quantum groupoids give a description of  
arbitrary finite index and finite depth II$_1$ subfactors via
a Galois correspondence and  
explain how to express subfactor invariants such as bimodule
categories and principal graphs in quantum groupoid terms. 
Recall that the {\em depth} (\cite{GHJ}, 4.6.4) 
of a subfactor $N\subset M$ with finite index $\lambda^{-1}=[M:N]$ is
\begin{eqnarray*}
n &=& \min\{ k\in \mathbb{Z}^+ \mid \dim Z(N'\cap M_{k-2}) = \dim
Z(N'\cap M_k) \} \in \mathbb{Z}_+ \cup \{\infty\},
\end{eqnarray*}
where $N\subset M\subset M_1\subset M_2\subset\cdots$
is the corresponding Jones tower.
 
\subsection{A Galois correspondence}
The following simple observation gives a natural passage from 
arbitrary finite depth subfactors to depth $2$ subfactors.

\begin{proposition}
\label{reducing depth}
For all $k\geq 0$ the inclusion $N\subset M_k$ has depth
$d+1$, where $d$ is the smallest positive integer $\geq
\frac{n-1}{k+1}$.
In particular, $N\subset M_i$ has depth $2$ for all $i\geq n-2$,
i.e., any finite depth subfactor is an intermediate
subfactor of some depth $2$ inclusion.
\end{proposition}
\begin{proof}
Note that $\dim Z(N'\cap M_i) = \dim Z(N'\cap M_{i+2})$ for all
$i\geq n-2$.
By \cite{PP2}, the tower of basic construction for $N\subset M_k$
is
$$
N\subset M_k \subset M_{2k+1} \subset M_{3k+2} \subset \cdots,
$$
therefore, the depth of this inclusion is equal to $d+1$, where $d$
is the smallest positive integer such that $d(k+1)-1 \geq n-2$.
\end{proof}

This result means that $N\subset M$ can be realized as an intermediate
subfactor of a crossed product inclusion $N\subset N\rtimes B$ for some
quantum groupoid $B$ :
$$
N\subset M \subset N\rtimes B.
$$
Recall that in the case of a Kac algebra
action there is a Galois correspondence between intermediate von Neumann
subalgebras of $N\subset N\rtimes B$ and left coideal $*$-subalgebras of
$B$ \cite{ILP},\cite{E1}. Thus, it is natural to ask about
a quantum groupoid analogue of this correspondence.

A unital $*$-subalgebra $I\subset B$ such that
$\Delta(I) \subset B\otimes I$
(resp.\ $\Delta(I) \subset I\otimes B$) is said to be
a left (resp.\ right) {\em coideal
$*$-subalgebra}.
The set $\ell(B)$ of left
coideal $*$-subalgebras is a lattice under the
usual operations:
$$
I_1 \wedge I_2 = I_1\cap I_2, \qquad I_1 \vee I_2 = (I_1\cup I_2)''
$$
for all $I_1,I_2\in \ell(B)$. The smallest element of $\ell(B)$ is $B_t$
and the greatest element is $B$. $I\in \ell(B)$ is said to be
{\em connected} if $ Z(I)\cap B_s =\mathbb{C}$.

To justify this definition, note that if $I=B$, then this is
precisely the definition of a connected quantum groupoid, 
and if $I=B_t$, then this definition is equivalent 
to $\widehat B$ being connected (\cite{N1}, 3.10, 3.11, \cite{BNSz}, 2.4).

On the other hand, the set 
$\ell(M_1\subset M_2)$ of intermediate von Neumann
subalgebras of $M_1\subset M_2$ also forms a lattice under the
operations
$$
K_1 \wedge K_2 = K_1\cap K_2, \qquad K_1 \vee K_2 = (K_1\cup K_2)''
$$
for all $K_1,K_2\in \ell(M_1\subset M_2)$. The smallest element
of this lattice is $M_1$ and the greatest element is $M_2$.

Given a left (resp.\ right)  action of $B$ on a von Neumann algebra
$N$, we will denote (by an abuse of notation)
$$
N\rtimes I = \span\{ [x\otimes b] \mid x\in N,\, b\in I\} \subset
N\rtimes B.
$$
The next theorem establishes a {\em Galois correspondence},
i.e., a lattice isomorphism between $\ell(M_1\subset M_2)$ and $\ell(B)$.

\begin{theorem}
\label{galois correspondence}
Let $N\subset M\subset M_1\subset M_2\subset\cdots$ be the tower
constructed from a depth $2$ subfactor $N\subset M$, $B=M'\cap M_2$
be the corresponding quantum groupoid, and $\theta$ be
the isomorphism between $M_1\rtimes B$ and $M_2$ 
(Proposition~\ref{crossproduct}).
Then
\begin{eqnarray*}
\phi &:& \ell(M_1\subset M_2) \to \ell(B) :
     K\mapsto \theta^{-1}(M'\cap K) \subset B\\
\psi &:& \ell(B) \to \ell(M_1\subset M_2) :
     I \mapsto \theta(M_1\rtimes I)\subset M_2.
\end{eqnarray*}
define isomorphisms between $\ell(M_1\subset M_2)$ and $\ell(B)$
inverse to each other.
\end{theorem}
\begin{proof}
First, let us check that $\phi$ and $\psi$ are indeed maps
between the specified lattices. It follows from the
definition of the crossed product that $M_1\rtimes I$ is a von
Neumann subalgebra of $M_1\rtimes B$, therefore $\theta(M_1\rtimes I)$
is a von Neumann subalgebra of $M_2=\theta(M_1\rtimes B)$,
so $\psi$ is a map to $\ell(M_1\subset M_2)$.
To show that $\phi$ maps to $\ell(B)$,
let us show that the annihilator $(M'\cap K)^0 \subset \widehat B$
is a left ideal in $\widehat B$.

For all $x\in A,\, y\in (M'\cap K)^0$, and $b\in M'\cap K$ we have
\begin{eqnarray*}
\la xy,\, b\ra
&=& \lambda^{-2}\tau(xye_2e_1wb) = \lambda^{-2}\tau(ye_2e_1wbx)\\
&=& \lambda^{-3} \tau(ye_2e_1E_{M'}(e_1wbx))
  = \la y,\,\lambda^{-1}E_{M'}(e_1wbx)\ra,
\end{eqnarray*}
and it remains to show that $E_{M'}(e_1wbx) \in M'\cap K$.
By (\cite{GHJ}, 4.2.7), the square
$$
\begin{array}{ccc}
K & \subset & M_2 \\
\cup &      & \cup \\
M'\cap K &\subset & M'\cap M_2
\end{array}
$$
is commuting, so $E_{M'}(K) \subset M'\cap K$.
Since $e_1wbx\in K$, then $xy \in (M'\cap K)^0$, i.e.,
$(M'\cap K)^0$ is a left ideal and $\phi(K)= \theta^{-1}(M'\cap K)$
is a left coideal $*$-subalgebra.

Clearly, $\phi$ and $\psi$ preserve $\wedge$ and $\vee$,
moreover $\phi(M_1) =B_t,\,\phi(M_2)=B$ and $\psi(B_t)=M_1,\,
\psi(B) = M_2$, therefore they are morphisms of lattices.

To see that they are inverses for each other, we first observe that
the condition $\psi\circ\phi =\id$ is equivalent to $M_1(M'\cap K)= K$,
and the latter follows from applying the conditional expectation $E_K$
to $M_1(M'\cap M_2)= M_1B = M_2$. The condition $\phi\circ\psi =\id$
translates into $\theta(I) = M'\cap \theta(M_1 \rtimes I)$.
If $b\in I,\, x\in M=M_1^B$, then
\begin{eqnarray*}
\theta(b)x
&=& \theta([1\otimes b][x\otimes 1])
     = \theta([(b\1\lact x) \otimes b\2]) \\
&=& \theta([x(1\2\lact 1)\otimes \eps(1\1 b\1)b\2 ])
     = \theta([x\otimes b]) =x\theta(b),
\end{eqnarray*}
i.e., $\theta(I)$ commutes with $M$. Conversely, if
$x\in M'\cap \theta(M_1 \rtimes I)\subset B$, then $x=\theta(y)$ for some
$y\in (M_1 \rtimes I) \cap B = I$, therefore $x\in \theta(I)$.
\end{proof}

\begin{corollary}
(\cite{NV3}, 4.5)
\begin{enumerate}
\item[(i)]
$K = M\rtimes I$ is a factor iff $ Z(I)\cap B_s =\mathbb{C}$.
\item[(ii)]
The inclusion $M_1 \subset K = M_1\rtimes I$ is irreducible iff
$B_s \cap I=\mathbb{C}$.
\end{enumerate}
\end{corollary}
\subsection{Bimodule categories}
Here we establish an equivalence between the tensor category
$\Bimod_{N-N}(N\subset M)$
of $N-N$ bimodules of a subfactor $N\subset M$ (which is, 
by definition, the tensor category generated by simple subobjects 
of ${}_N L_2(M_n)_N,\, n\geq -1$, where $M_0=M$  and $M_{-1}=N$)
and the co-representation category of the $C^*$-quantum groupoid $B$
canonically  associated with it as in Theorem~\ref{galois correspondence}.

First introduce several useful categories associated to $B$.
Recall that a left (resp., right) $B$-comodule $V$ (with the structure map
denoted by $v \mapsto v\I \otimes v\II,\, v\in V$)
is said to be {\em unitary}, if
\begin{eqnarray*}
(v_2\I)^*(v_1, v_2\II) &=& S(v_1\I)g (v_1\II,v_2) \\
(\mbox{resp., } (v_1\II)(v_1\I, v_2) &=& g^{-1}S((v_1\I)^*)(v_1,v_2\I)),
\end{eqnarray*}
where $v_1,v_2\in V,$ and $g$ is the canonical group-like element of
$B$. This definition for Hopf
$*$-algebra case can be found, e.g., in (\cite{KS}, 1.3.2).

Given left coideal $*$-subalgebras $H$ and $K$ of $B$, we consider
a category $\C_{H-K}$ of left relative $(B, H-K)$ Hopf bimodules
(cf. \cite{Ta}), whose objects are Hilbert spaces which
are both $H-K$-bimodules and left unitary $B$-comodules such that the bimodule
action commutes with the coaction of $B$, i.e., for
any object $V$ of $\C_{H-K}$ and $v\in V$ one has
$$
(h\lact v \ract k)\I \otimes (h\lact v \ract k)\II =
h\1 v\I k\1 \otimes (h\2 \lact v\II \ract k\2),
$$
where $v\mapsto v\I \otimes v\II$ denotes the coaction of $B$ on
$v$, $h\in H,\ k\in K$, and morphisms are  intertwining maps.

Similarly one can define a category of right relative
$(B, H-K)$ Hopf bimodules.

\begin{remark}
\label{on relative bimodules}
Any left $B$-comodule $V$ is automatically a $B_t-B_t$-bimodule via
$z_1\cdot v\cdot z_2 =  \eps(z_1v\I z_2)v\II,~v\in V, z_1,z_2\in B_t$.
For any object of $\C_{H-K}$, this $B_t-B_t$-bimodule structure
is a restriction of the given $H-K$-bimodule structure; it is easily seen
by applying $(\eps\otimes id)$ to both sides of the above 
relation of commutation and taking $h,k\in B_t$.
Thus, if $H=B_t$ (resp.\ $K=B_t$),
we can speak about right (resp. \ left) relative Hopf modules,
a special case of weak Doi-Hopf modules \cite{Bo}.
\end{remark}
\begin{proposition}
\label{double cosets}
If $B$ is a  group Hopf $C^*$-algebra and $H,K$ are subgroups,
then there is a bijection between simple objects of $\C_{H-K}$
and double cosets of $H\backslash B/K$.
\end{proposition}
\begin{proof}
If $V$ is an object of $\C_{H-K}$, then every simple subcomodule
of $V$ is $1$-dimensional. Let $U =\mathbb{C}u$ 
($u\mapsto g\otimes u,\, g\in B$)
be one of these comodules, then all other simple subcomodules of $V$
are of the form $h\lact U \ract k$, where $h\in H,\,k\in K$, and
$$
V= \oplus_{h,k}\, (h\lact U \ract k) = \span\{ HgK \}.
$$
Vice versa, $\span\{ HgK \}$ with natural $H-K$ bimodule and $B$-comodule
structures is a simple object of $\C_{H-K}$.
\end{proof}
\begin{example}
\label{basic object}
If $H,V,K$ are left coideal $*$-subalgebras of $B,\ H\subset V,K\subset V$,
then $V$ is an object of $\C_{H-K}$ with the structure maps given by
$ h\lact v \ract k = hvk$ and $\Delta$, where
$v\I \otimes v\II =\Delta(v),\ v\in V,\ h\in H,\ k\in K$. The scalar
product is defined by the restriction on $V$ of the {\em Markov trace} of 
the connected inclusion $B_t\subset B$ (\cite{JS}, 3.2).
Similarly, right coideal $*$-subalgebras of $B$ give examples of right relative
$(B, H-K)$ Hopf bimodules.
\end{example}

Given an object $V$ of $\C_{H-K}$, it is straightforward to show that
the conjugate Hilbert space
$\overline{V}$ is an object of $\C_{K-H}$ with the bimodule action
$$
k\lact\overline{v}\ract h =\overline{h^*\lact v \ract k^*}\ \ \ 
(\forall h\in H,\, k\in K)
$$
(here $\overline{v}$ denotes the vector $v\in V$ considered as an
element of $\overline{V}$) and the coaction $\overline{v}\mapsto \overline{v}\I
\otimes \overline{v}\II=(v^{(1)})^*\otimes \overline{v}\II$. 
Define $V^*$, the dual object of $V$, to be $\overline{V}$ with
the above structures.
One can directly check that $V^{**}\cong V$ for any object  $V$.
In Example \ref{basic object} the dual object can
be obtained by putting $\overline{v}=v^*$ for all $v\in V$.
\medskip
\begin{definition}
\label{reltenspr}
Let $L$ be another coideal $*$-subalgebra of $B$. For any objects
$V\in \C_{H-L}$ and $W\in \C_{L-K}$, we define an object $V\otimes_L W$
from $\C_{H-K}$ as a tensor product of bimodules $V$ and $W$ \cite{JS}, 4.1
equipped with a comodule structure
$$
(v\otimes_L w)^{(1)}\otimes (v\otimes_L w)^{(2)}=v^{(1)}w^{(1)}\otimes
(v^{(2)}\otimes_L w^{(2)}).
$$
\end{definition}
One can verify that we have indeed an object 
from $\C_{H-K}$ and that
the operation of tensor product is
(i)
associative, i.e., $V\otimes_L (W \otimes_P U)\cong(V\otimes_L W)\otimes_P U$;
(ii)
compatible with duality, i.e., $(V\otimes_L W)^* \cong W^*\otimes_L V^*$;
(iii)
distributive, i.e.,
$(V\oplus V')\otimes_L W = (V\otimes_L W) \oplus (V'\otimes_L W)$
(\cite {NV3}, 5.5).
The tensor product of morphisms $T\in \Hom(V, V')$ and $S\in
\Hom(W, W')$ is defined as usual :
$$
(T\otimes_L S)(v\otimes_L w) = T(v)\otimes_L S(w).
$$

From now on let us suppose that $B$ is biconnected and
acts outerly on the left on a II${}_1$ factor $N$.
Given an object $V$ of $\C_{H-K}$, we construct an
$N\rtimes H-N\rtimes K$-bimodule $\widehat{V}$ as follows.  We put
$$
\widehat{V} = \span \{ \Delta(1)\lact (\xi\otimes v)\mid
\xi\otimes v\in L^2(N)\otimes_\mathbb{C} V\}= L^2(N)\otimes_{B_t} V
$$
and  denote $[\xi\otimes v] = \Delta(1)\lact(\xi\otimes v)$.
Let us equip $\widehat{V}$ with the scalar product
$$
([\xi\otimes v],[\eta\otimes w])_{\widehat{V}}
=(\xi,\eta)_{\widetilde{L^2(N)}}(v,w)_V.
$$
and define the actions of $N,H,K$ on $\widehat{V}$ by
$$
a[\xi\otimes v] = [a\xi\otimes v], \qquad
[\xi\otimes v]a = [\xi(v\I\lact a)\otimes v\II],
$$
$$
h[\xi\otimes v] = [(h\1\lact \xi)\otimes (h\2\lact v)], \qquad
[\xi\otimes v]k = [\xi\otimes (v\ract k)],
$$
for all $a\in N,\, h\in H,\, k\in K$. One can check that these
actions define the structure of an
$(N\rtimes H)-(N\rtimes K)$ bimodule on $\widehat{V}$ in
the algebraic sense and that this bimodule is unitary.

For any morphism $T\in \Hom(V,W)$, define 
$\widehat{T} \in \Hom(\widehat{V},\widehat{W})$ by
$$
\widehat{T}([\xi\otimes v]) = [\xi\otimes T(v)].
$$
\begin{example}
\label{crossproduct as bimodule}
(\cite{NV3}, 5.6)
For $V\in \C_{H-K}$ from Example \ref{basic object}, we have
$\widehat{V}=N\rtimes V$ as $N\rtimes H-N\rtimes K$-bimodules.
\end{example}
The proof of the following theorem is purely technical 
(see \cite{NV3}, 5.7):
\begin{theorem}
\label{properties of the functor}
The above assignments $V\mapsto \widehat{V}$ and $T\mapsto \widehat{T}$
define a functor from  $\C_{H-K}$ to the category of
$N\rtimes H - N\rtimes K$ bimodules. This functor preserves direct
sums and is compatible with operations of taking tensor products and
adjoints in the sense that if $W$ is an object of $\C_{K-L}$, then
$$
\widehat {V\otimes_K W} \cong \widehat{V} \otimes_{N\smrtimes K} \widehat{W},
\quad \mbox{ and } \quad \widehat{V^*} \cong (\widehat{V})^*.
$$
\end{theorem}
According to Remark~\ref{on relative bimodules},  $\C_{B_t-B_t}$
is nothing but the category of $B$-comodules, $\Corep(B)$.
Let us show that this category is equivalent to $\Bimod_{N-N}(N\subset M)$.
\begin{theorem}
\label{N-N equivalence}
Let $N\subset M$ be a finite depth subfactor with finite index, $k$
be a number such that $N\subset M_k$ has depth $\leq 2$,
and let $B$ be a canonical quantum groupoid
such that $(N\subset M_k)\cong (N\subset N\rtimes B)$.
Then  $\Bimod_{N-N}(N\subset M)$ and  $\Rep(B^*)$
are equivalent as tensor categories.
\end{theorem}
\begin{proof}
First, we observe that
$$
\Bimod_{N-N}(N\subset M) = \Bimod_{N-N}(N\subset M_l)
$$
for any $l\geq 0$. Indeed, since both categories are semisimple, it is
enough to check that they have the same set of simple objects.
All objects of $\Bimod_{N-N}(N\subset M_l)$ are also objects of
$\Bimod_{N-N}(N\subset M)$. Conversely, since irreducible $N-N$ subbimodules
of ${}_N L^2(M_i)_N$ are contained in the decomposition of
${}_N L^2(M_{i+1})_N$ for all $i\geq 0$, we see that objects of
$\Bimod_{N-N}(N\subset M)$ belong to $\Bimod_{N-N}(N\subset M_l)$.

Hence, by Proposition~\ref{reducing depth},
it suffices to consider
the problem in the case when $N\subset M$ has
depth $2$ ($M=N\rtimes B$), i.e., to prove that
$\Bimod_{N-N}(N\subset N\rtimes B)$ is equivalent to $\Corep(B)$.
Theorem~\ref{properties of the functor} gives a functor from
$\Corep(B) = \Rep(B^*)$ to $\Bimod_{N-N}(N\subset N\rtimes B)$.
To prove that this functor is an equivalence,
let us check that it yields a bijection
between classes of simple objects of these categories.

Observe that $B$ itself is an object of $\Corep(B)$ via
$\Delta : B\to B\otimes_\mathbb{C} B$ and $\widehat{B} = {}_N L^2(M)_N$.
Since the inclusion $N\subset M$ has depth $2$, the simple
objects of $\Bimod_{N-N}(N\subset M)$
are precisely irreducible subbimodules of ${}_N L^2(M)_N$.
We have $\widehat{B} = {}_N L^2(M)_N = \oplus_i\,   {}_N
(p_iL^2(M))_N$, where $\{p_i\}$ is a family of mutually orthogonal
minimal projections in $N'\cap M_1$ so every
bimodule $p_iL^2(M)$ is irreducible. On the other hand, $B$ is
cosemisimple, hence $B=\oplus_i\, V_i$, where each $V_i$
is an irreducible subcomodule. Note that $N'\cap M_1= B^* =\sum p_iB^*$
and every $p_iB$ is a simple submodule of $B$ ($=$ simple subcomodule
of $B^*$). Thus, there is a bijection between the sets of simple
objects of $\Corep(B)$ and $\Bimod_{N-N}(N\subset M)$, so
the categories are equivalent.
\end{proof}

\subsection{Principal graphs}
The principal graph of a subfactor $N\subset M$
is defined as follows (\cite{JS}, 4.2, \cite{GHJ}, 4.1).
Let $X= {}_N L^2(M)_M$ and consider the following sequence of
$N-N$ and $N-M$ bimodules :
$$
{}_N L^2(N)_N,~X,~X\otimes_M X^*,~X\otimes_M X^*\otimes_N X,\dots
$$
obtained by right tensoring with $X^*$ and $X$. The vertex set
of the principal graph is indexed by the classes of simple bimodules
appearing as summands in the above  sequence. Let us connect
vertices corresponding to bimodules ${}_N Y_N$ and ${}_N Z_M$ by $l$
edges if ${}_N Y_N$ is contained in the decomposition of ${}_N Z_N$,
the restriction of ${}_N Z_M$, with multiplicity $l$.


We apply Theorem~\ref{properties of the functor} to express the principal
graph of a finite depth subfactor $N\subset M$ in terms of the
quantum groupoid associated with it.

Let $B$ and $K$ be a quantum groupoid and its left coideal
$*$-subalgebra such that $B$ acts on $N$ and $(N\subset M) \cong
(N\subset N\rtimes K)$. Then $\widehat{K} = {}_N L^2(M)_M$, where
we view $K$ as a relative $(B, K)$ Hopf module as in
Example~\ref{basic object}

By Theorem~\ref{properties of the functor} we can identify
irreducible $N-N$ (resp.\ $N-M$) bimodules with  simple $B$-comodules
(resp.\  relative right $(B, K)$ Hopf modules). Consider a bipartite
graph with vertex set given by the union of (classes of) simple
$B$-comodules and simple relative right $(B, K)$ Hopf modules
and the number of edges between the vertices $U$ and $V$ representing
$B$-comodule and relative right $(B, K)$ Hopf module respectively
being  equal to the multiplicity of $U$ in the decomposition of $V$
(when the latter is viewed as a $B$-comodule) :
\begin{eqnarray*}
& \mbox{simple}~B\mbox{-comodules}& \\
& \cdots \vert \cdots \vert \cdots & \\
& \mbox{simple relative right}~(B,K)~\mbox{Hopf modules} &
\end{eqnarray*}
The principal graph of $N\subset M$ is the connected part of the above graph
containing the trivial $B$-comodule.


It was shown in (\cite{NV3}, 3.3, 4.10) that:
\begin{enumerate}
\item[(a)] the map $\delta : I \mapsto [g^{-1/2}S(I)g^{1/2}]'\cap \widehat B
\subset \widehat B\rtimes B$ defines a lattice anti-isomorphism between  
$\ell(B)$ and $\ell(\widehat B)$;
\item[(b)] the triple $\delta(I) \subset \widehat B \subset \widehat B\rtimes I$
is a basic construction.
\end{enumerate}

\begin{proposition}
\label{Bratteli diagram = principal graph}
If $K$ is a coideal $*$-subalgebra of $B$ then the principal
graph of the subfactor $N\subset N\rtimes K$ is given by
the connected component of the Bratteli
diagram of the inclusion $\delta(K) \subset \widehat B$
containing the trivial representation of $\widehat B$.
\end{proposition}
\begin{proof}
First, let us show that there is a bijective correspondence
between right relative $(B,K)$ Hopf modules and $(\widehat B\rtimes K)$-modules.
Indeed, every right $(B,K)$ Hopf module $V$ carries a right action of $K$.
If we define a right action of $\widehat B$ by
$$
v\ract x = \la v\I,\, x\ra v\II,\quad v\in V, x\in \widehat B,
$$
then we have
\begin{eqnarray*}
(v \ract k)\ract x
&=& \la v\I k\1,\,x \ra (v\II \ract k\2) 
 = \la v\I,\, (k\1\lact x) \ra (v\II \ract k\2) \\
&=& (v \ract (k\1\lact x)) \ract k\2,
\end{eqnarray*}
for all $x\in \widehat B$ and $k\in K$ which shows that $kx$
and $(k\1\lact x)k\2$ act on $V$ exactly in the same way,
therefore $V$ is a right $(\widehat B\rtimes K)$-module.

Conversely, given an action of $(\widehat B\rtimes K)$ on $V$,
we automatically have a $B$-comodule structure such that
\begin{eqnarray*}
\la v\I k\1,x \ra (v\II \lact k\2)
&=& \la v\I,\, x\1\ra \la k\1,\, x\2 \ra (v\II \ract k\2) \\
&=& (v\ract (k\1\lact x)) \ract k\2 = (v\ract k)\ract x \\
&=& \la x,\, (v \ract k)\I \ra (v \ract k)\II,
\end{eqnarray*}
which shows that $v\I k\1 \otimes (v\II \ract k\2) =
(v \ract k)\I \otimes (v \ract k)\II$, i.e., that $V$
is a right relative $(B,K)$-module.

Thus, we see that the principal graph is given by
the connected component the Bratteli diagram
of the inclusion $\widehat B \subset \widehat B\rtimes K$
containing the trivial representation of $\widehat B$.
Since $\widehat B\rtimes K$
is the basic construction for the inclusion $\delta(K)\subset \widehat B$,
therefore the Bratteli diagrams of the above two inclusions are the same.
\end{proof}

\begin{corollary}
\label{principal graph for depth 2}
If $N\subset N\rtimes B$ is a depth 2 inclusion
corresponding to the quantum groupoid $B$, then
its principal graph is given by the Bratteli diagram
of the inclusion $\widehat B_t \subset \widehat B$.
\end{corollary}
\begin{proof}
In this case $K=B$ and inclusion $\widehat B_t \subset \widehat B$
is connected, so that $\delta(K) = \widehat B_t$ (note that
$\widehat B$ is biconnected).
\end{proof}
\end{section}

\bibliographystyle{amsalpha}

\end{document}